\documentclass[12pt]{amsart}
\usepackage{amsmath,amsfonts,amsthm,amscd,amssymb,mathrsfs,amssymb}
\usepackage{graphicx}
\newtheorem{theorem}{Theorem}

\newtheorem{proposition}[theorem]{Proposition}
\newtheorem{lemma}[theorem]{Lemma}
\newtheorem{definition}[theorem]{Definition}

\newtheorem{corollary}[theorem]{Corollary}

\newtheorem{remark}[theorem]{Remark}

\renewcommand{\H}{\mathcal{H}}
\newcommand{\CP}{\mathbb{CP}}
\newcommand{\CC}{\mathbb{C}}

\newcommand{\RR}{\mathbb{R}}
\newcommand{\ZZ}{\mathbb{Z}}

\newcommand{\U}{{\rm{U}}}
\newcommand{\LB}{{\rm{LB}}}
\newcommand{\GH}{{\rm{GH}}}
\renewcommand{\i}{i}
\newcommand{\w}{\omega}

\numberwithin{equation}{section}
\numberwithin{theorem}{section}
\begin{document}
\bibliographystyle{amsalpha} 
\title[Monopole metrics]{Monopole metrics and the orbifold
\\[2pt]Yamabe problem} 
\author{Jeff A. Viaclovsky}
\address{Department of Mathematics, University of Wisconsin, Madison, 
WI, 53706}
\email{jeffv@math.wisc.edu}
\dedicatory{To Pierre B\'erard and Sylvestre Gallot on the 
occasion of their sixtieth birthdays.}
\thanks{Research partially supported by NSF Grant DMS-0804042}
\begin{abstract}
We consider the self-dual conformal classes 
on $n \# \CP^2$ discovered by LeBrun. These depend upon 
a choice of $n$ points in hyperbolic $3$-space, called monopole 
points. We investigate the limiting behavior of various constant scalar 
curvature metrics in these conformal classes as the 
points approach each other, or as the points tend to the boundary 
of hyperbolic space. 
There is a close connection to the orbifold Yamabe problem,
which we show is not always solvable (in contrast to the 
case of compact manifolds). In particular, we show that there 
is no constant scalar curvature orbifold metric in the conformal 
class of a conformally compactified non-flat hyperk\"ahler ALE space
in dimension four. 
\end{abstract}
\date{February 10, 2010. Revised October 2010.}
\maketitle
\parskip1pt
\setcounter{tocdepth}{1}
\vspace{-5mm}
\tableofcontents

\section{Introduction}

There is an interesting history regarding the existence 
of self-dual metrics on $n \# \CP^2$ beginning with 
work of Yat-Sun Poon \cite{Poon}.
Using techniques from twistor theory, Poon proved the existence 
of a $1$-parameter family of self-dual conformal classes on 
$\CP^2 \# \CP^2$ and that any such conformal class 
with positive scalar curvature must be in this family. 
Examples for larger $n$ were found by  Donaldson-Friedman \cite{DonaldsonFriedman}
and Floer \cite{Floer} using gluing methods. 
In 1991, Claude LeBrun \cite{LeBrunJDG} produced explicit examples 
with ${\rm{U}}(1)$-symmetry on $n \# \CP^2$, using a hyperbolic ansatz 
inspired by the Gibbons-Hawking ansatz \cite{GibbonsHawking}.   
LeBrun's construction depends on the choice of $n$ points 
in hyperbolic 3-space $\mathcal{H}^3$. For $n = 2$, the only 
invariant of the configuration is the distance 
between the monopole points, and LeBrun conformal classes 
are the same as the $1$-parameter family found by Poon. 

\subsection{Limits of LeBrun metrics}
The first question we address in this paper: is there a nice 
compactification of the moduli space of LeBrun metrics 
on  $n \# \CP^2$?
In general, as the monopole points limit towards each other, 
or if the points approach the boundary of hyperbolic space, 
some degeneration will occur. 
We emphasize that the LeBrun construction produces {\em{conformal classes}}
on $n \# \CP^2$. To discuss convergence in the Cheeger-Gromov sense, 
one needs to choose a conformal factor.
Of course, the limit will strongly 
depend on the particular choice of conformal metrics. 
Some degenerations were already described in \cite{LeBrunJDG}
and \cite{DonaldsonFriedman}, but these examples depended on a somewhat 
arbitrary choice of conformal factor. 

The solution of the Yamabe problem 
provides one with a very natural metric in these conformal classes. 
However, the abstract existence theorem does not tell one what the actual 
minimizer looks like in any particular case, and 
other methods are needed to understand the geometry 
of minimizers. The main 
point of this paper is to describe the limiting behavior
of the {\em{Yamabe minimizers}} in these conformal classes
as they degenerate. In general, Yamabe minimizers are not necessarily unique; 
an example of non-uniqueness is given Theorem \ref{n=2thm}.
We also examine the existence and limiting 
behavior of various non-minimizing constant scalar curvature metrics. 
Given a subgroup of the conformal
automorphism group, Hebey-Vaugon have shown there is a minimizer of the 
Yamabe functional when restricted to the class of 
invariant functions (the equivariant Yamabe problem), 
and these automorphisms will act as {\em{isometries}} on 
the minimizer \cite{Hebey1996}. Of course, a symmetric Yamabe minimizer 
can have higher energy than a Yamabe minimizer, and an 
example of this is seen in Theorem \ref{n=2thm}.   

If $G \subset {\rm{SO}}(4)$ is a finite subgroup acting freely 
on $S^3$, then we let
$G$ act on $S^4 \subset \RR^5$ acting as rotations 
around the $x_5$-axis. The quotient $S^4/G$ is then a
orbifold, with two singular points, and the 
spherical metric $g_{S}$ descends to this orbifold. 
Near the singular points, 
the metric is asymptotic to a cone metric $\mathcal{C}( S^3/G)$,
thus $S^4/G$ looks like a United States ``football''. 
In the following, $G \subset {\rm{SU}}(2)$ will be 
a certain cyclic subgroup $\ZZ_m$, see \eqref{su2} below. 
For a smooth Riemannian manifold $(M,g)$, the Yamabe invariant of the 
conformal class is denoted by $Y(M, [g])$, see Section~\ref{YI}.
If $(M,g)$ is an orbifold, then $Y_{orb}(M, [g])$ will denote the orbifold 
Yamabe invariant, see Section~\ref{OYP}. 

\begin{figure}
\includegraphics{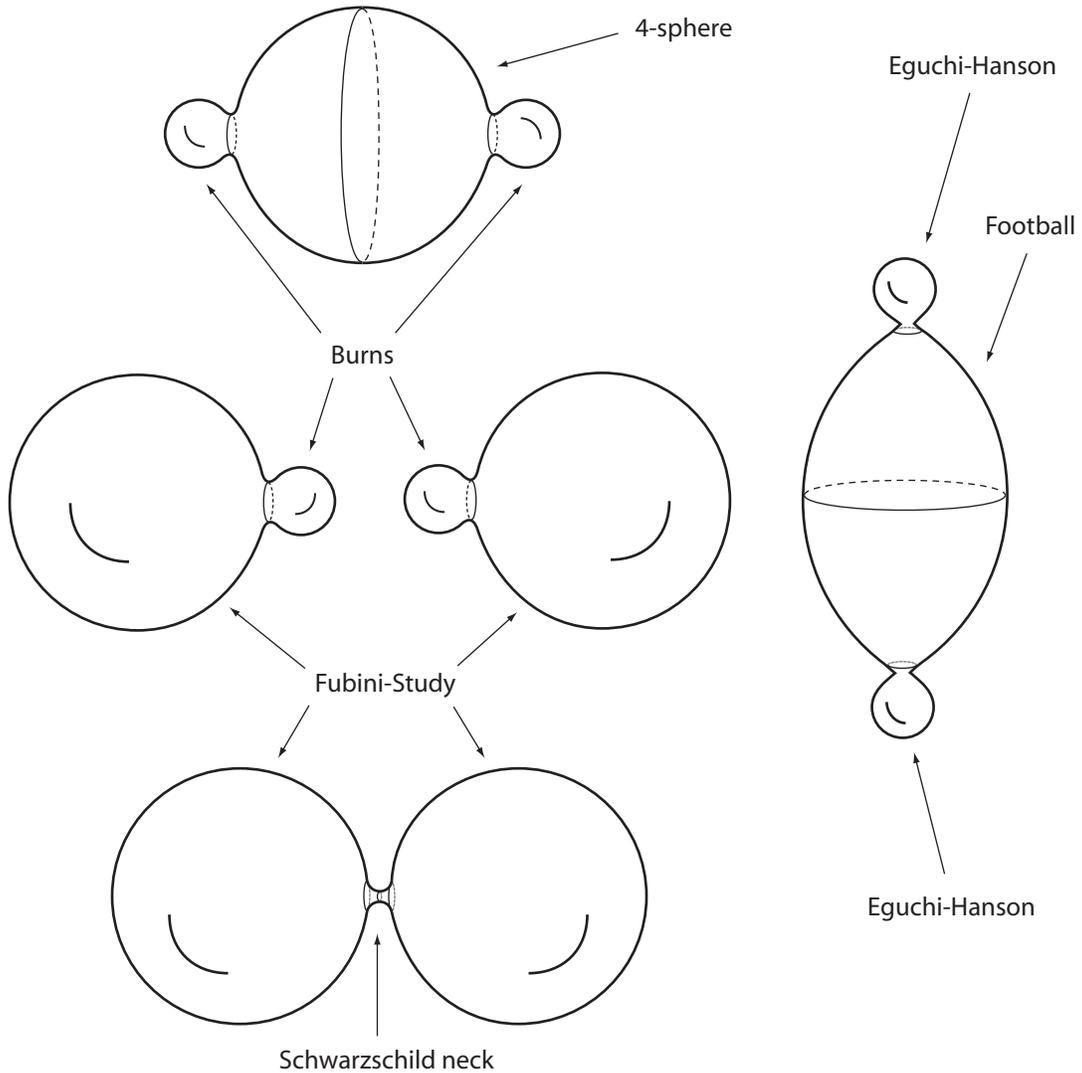}
\caption{The 5 limiting cases described in Theorem \ref{n=2thm}.
The left side is as $d(p_1,p_2) \rightarrow \infty$.  
The top left is the symmetric $K$-Yamabe minimizer. The middle 
left are the two Yamabe minimizers. The lower left 
is the metric obtained by Joyce gluing. The right hand side 
is as $d(p_1,p_2) \rightarrow 0$. In this case, the Yamabe 
minimizers and the symmetric Yamabe minimizers have the same limiting behavior.
Note that in each of the above cases, to obtain $\CP^2 \# \CP^2$ topologically,
the attaching map for one of the factors should be orientation reversing.} 
\label{bubble2}
\end{figure}

We first discuss the special case of $n=2$. 
To employ the equivariant Yamabe problem,  
one must first understand the group of conformal 
automorphisms: 
it was proved in \cite{HV} that the conformal 
group $G$ of Poon's metrics for $n=2$ is given by
\begin{align}
G = ({\rm{U}}(1) \times {\rm{U}}(1)) \ltimes D_4, 
\end{align}
where $D_4$ is the dihedral group of order $8$. 
There is the index $2$ subgroup given by 
\begin{align}
K = ({\rm{U}}(1) \times {\rm{U}}(1)) \ltimes (\mathbb{Z}_2 \oplus \mathbb{Z}_2), 
\end{align}
which are exactly the lifts of hyperbolic isometries preserving 
the set of $2$ monopole points. 
In contrast, for $n > 2$, any conformal automorphism 
of a LeBrun metric is a lift of an 
isometry of $\mathcal{H}^3$. There is an ``extra'' involution when 
$n =2$, which is not a lift of any hyperbolic isometry, see \cite{HV}.
Let $d_H(\cdot, \cdot)$ denote hyperbolic distance. 
\begin{theorem}
\label{n=2thm}
Let $(M,g)$ be a Poon-LeBrun metric on $\CP^2 \# \CP^2$
with monopole points $p_1$ and $p_2$. 
The Yamabe invariant satisfies the sharp estimate
\begin{align}
\label{2est}
8 \pi \sqrt{3} = Y_{orb}(S^4/ \ZZ_2, [g_S]) < Y(M,[g]) 
<  Y(\CP^2, [g_{FS}]) = 12 \pi \sqrt{2}. 
\end{align}
There exists a number $N$ large, such that if
$d_H(p_1, p_2) > N$ then the following holds. 
There are two distinct Yamabe minimizers, 
each limiting to $g_{FS}$ on $\CP^2$ as $d_H(p_1, p_2) \rightarrow \infty$.  
In each case, there is one singular point of convergence at which a 
Burns metric bubbles off. 
The symmetric $K$-Yamabe minimizers limit 
to $(S^4, g_{S})$  as  $d_H(p_1, p_2) \rightarrow \infty$,
with 2 antipodal singular points of 
convergence, with Burns metrics bubbling off 
at each of the singular points. 
There is a fourth constant 
scalar curvature metric, limiting to  $\CP^2 \vee \CP^2$
(the wedge of two copies of $\CP^2$ with Fubini-Study metrics, 
touching at a single point) as $d_H(p_1, p_2) \rightarrow \infty$.
In this case, a Euclidean Schwarzschild metric with two asymptotically 
flat ends bubbles off.

  As $d_H(p_1, p_2) \rightarrow 0$, 
the limit of both the Yamabe minimizers and
the symmetric $K$-Yamabe minimizers is the $S^4 / \ZZ_2$-football 
with the round metric. In both cases, at each singular point an Eguchi-Hanson 
metric bubbles off. 
\end{theorem}
These limits are illustrated in Figure \ref{bubble2}.
An important point is that the properties of being self-dual and 
having constant scalar curvature form an elliptic system \cite{TV}.
In Tian-Viaclovsky \cite{TV2, TV3}, it was shown that 
limits of such metrics may have at worst multi-fold singularities, 
provided that the sequence has bounded $L^2$-norm of curvature 
and does not collapse. The crucial ingredient of this theory 
is the upper volume growth estimate proved in \cite{TV}. 
For other results dealing with this 
type of convergence in various settings, see 
\cite{AkutagawaDGA, AkutagawaPacific, Anderson, Andersonc, Bando, 
ChangQingYang, ChenWeber, ChenLeBrunWeber, Nakajimaconv, Tian}. 
In the situation considered in this paper, the $L^2$-curvature bound follows 
from the Chern-Gauss-Bonnet 
formula and Hirzebruch signature theorem. The non-collapsing condition 
will follow from uniform positivity of the Yamabe invariant.

 As mentioned above, the LeBrun monopole construction depends upon the choice 
of $n$ points in hyperbolic space. An easy generalization 
of this construction allows one to assign integer multiplicities 
greater than one at the monopole points. The resulting 
space will have orbifold points. We call such a 
space a {\em{LeBrun orbifold}}. 
Next, for $n > 2$, we present the following compactness 
theorem. 
\begin{theorem} 
\label{maint}
Let $(M,g)$ be a LeBrun self-dual conformal class
on $n \# \CP^2$ with monopole points $\{ p_1, \dots, p_n\}$. 
Then 
\begin{align}
\label{uest}
Y(M, [g]) \leq  Y(\CP^2, [g_{FS}]) = 12 \pi \sqrt{2},
\end{align}
with strict inequality for $n \geq 2$, 
where $g_{FS}$ denotes the Fubini-Study metric. 
Next, assume that all monopole points are contained in a 
compact set $\mathcal{K} \subset \mathcal{H}^3$. Then
there exists a constant $\delta_n$ depending 
only upon $n,\mathcal{K}$ such that
\begin{align}
\label{mest}
0 < \delta_n \leq Y(M, [g]).
\end{align}
Furthermore, any sequence of unit volume Yamabe minimizers in a sequence of 
LeBrun conformal classes (for fixed $n$) satisfying \eqref{mest} has a 
subsequence which converges (in the Cheeger-Gromov sense) 
to either to (1) a compactified LeBrun orbifold metric 
with $1 \leq k \leq n$ points or (2) the round metric on $S^4 / \ZZ_m$, 
for some $2 \leq m \leq n$. 
In addition, the estimate \eqref{mest} is true for $n = 2, 3$
without the requirement that the points are contained in a 
compact set $\mathcal{K}$. 
\end{theorem}

There can exist sequences of Yamabe minimizers 
with limiting behavior as in Case (1); this limit can
occur when some of the monopole points limit to the 
boundary of $\H^3$, as seen in Theorem \ref{n=2thm}. Another 
example is given in Theorem \ref{imaint} below.
There also can exist sequences limiting as in Case (2), as seen in Theorem \ref{n=2thm}. 

The estimate for the lower bound in \eqref{mest} is not explicit, 
this is proved by a contradiction argument in Section \ref{convergence}. 
It would be very interesting to 
find a sharp constant. We conjecture that 
$\delta_n = Y_{orb}( S^4 / \ZZ_n, [g_{S}])$,
without any requirement that the   
monopole points are contained in a compact set $\mathcal{K}$. 
For $n = 2,3$, the uniform positivity of the 
Yamabe invariant holds for topological reasons, see 
Proposition \ref{23prop}. 

We mention that the degeneration of the 
LeBrun conformal classes can also be studied using twistor 
theory. For this important perspective, we refer the reader to 
the recent paper of Nobuhiro Honda \cite{Hondadegeneration}. 
\subsection{The orbifold Yamabe problem}
Akutagawa and Botvinnik considered the Yamabe problem on 
orbifolds in \cite{AB1, AB2} in which they proved several 
foundational results. We will describe this in more detail 
in Section \ref{OYP}. Since the limits described above 
are typically orbifolds, it is no surprise that there is a 
close connection with the orbifold Yamabe problem. 
In fact, an important tool in identifying the possible limit spaces above 
is the following nonexistence result. 
\begin{theorem}
\label{nonexist}
Let $(X_n,g)$ be a hyperk\"ahler ALE metric in
dimension $4$, with group $G$ of order $n >1$ at infinity, 
and let $(\hat{X}, [\hat{g}])$ denote 
the orbifold conformal compactification. 
Then $Y_{orb}(\hat{X}, [\hat{g}]) = Y_{orb} ( S^4 / G, [g_S])$,
and there is no solution to the orbifold Yamabe
problem on $(\hat{X}, [\hat{g}])$. That is, there
is no conformal metric $\tilde{g} = e^{2u} \hat{g}$ 
having constant scalar curvature. 
\end{theorem}
This will be proved in Section \ref{OYP} along with some
other remarks on the orbifold Yamabe problem. 
The ALE metrics above, can be viewed as the Green's 
function metrics $g_p = \Gamma_p^2 \hat{g}$ of the orbifold 
compactification, where $\Gamma_p$ is the 
Green's function for the conformal Laplacian of $\hat{g}$ based 
at the orbifold point $p$. It is interesting that these  
ALE spaces have zero mass 
and their compactifications 
do {\em{not}} admit a solution of the orbifold Yamabe
problem. This shows that the orbifold Yamabe problem is 
more subtle than in the case of smooth manifolds.

In Section 7, we will prove another nonexistence result regarding
the negative mass ALE spaces found in \cite{LeBrunnegative}. 
\begin{theorem}
\label{nonsym}
If $(X,g)$ is a LeBrun negative mass ALE metric on $\mathcal{O}(-n)$
with $n > 1$,  then there is no symmetric solution
of the orbifold Yamabe problem on 
$(\hat{X}, \hat{g})$ invariant under $\rm{SU}(2)$. 
\end{theorem}
As a consequence, the symmetric Yamabe problem on orbifolds is 
not always solvable either. For $n =2$, this metric is 
the same as the compactified Eguchi-Hanson metric, which does not 
admit {\em{any}} constant scalar curvature metric by Theorem \ref{nonexist}.
We do not know if there is a non-symmetric solution on these orbifolds 
for $n \geq 3$. 

Finally, we present an existence result for the orbifold 
Yamabe problem, which for simplicity we state here in only the 
case of total multiplicity $3$ (see Corollary \ref{orbex}
for the general statement). The following theorem also shows 
that LeBrun {\em{orbifold}} metrics can in fact arise 
as a limit of smooth Yamabe metrics in LeBrun 
conformal classes, and also $k$ can be strictly less than 
$n$ in Case (1) of Theorem \ref{maint}. 
\begin{theorem}
\label{imaint}
Let $(M, [g])$ be compact self-dual LeBrun orbifold corresponding 
to a monopole point $p_1$ of multiplicity $1$ and a monpole 
point $p_2$ of multiplicity $2$, with $p_1 \neq p_2$.  
Then there exists a radius $r > 0$ such that if 
$d_{\H} (p_1, p_2) < r$, then  $(M, [g])$
admits a solution to the orbifold Yamabe problem.  
Furthermore, let $(M_1,[g_j])$ be a self-dual LeBrun conformal 
class on $3 \# \CP^2$ 
corresponding to $3$ distinct monopole points all of multiplicity one, 
with $p_1$ and $p_2$ fixed, and $p_{2}, p_{j,3} \in B(p_1, r)$.  
Then as $p_{j,3} \rightarrow p_{2}$, a subsequence of 
Yamabe minimizers on $(M_1,[g_j])$ converges to an orbifold 
Yamabe metric on $(M, [g])$. There is one singular point
of convergence, at which an Eguchi-Hanson metric bubbles off. 

 Next, let $p_1$ and $p_2$ be fixed and let $p_{j,3}$ limit to the 
boundary of $\H^3$ as $j \rightarrow \infty$.  
Then any sequence of Yamabe minimizers on $(M_1, [g_j])$ 
has a subsequence which converges to a Yamabe minimizer 
in the 2-pointed smooth LeBrun conformal class $(2 \# \CP^2,[\hat{g}_{LB} (p_1,p_2)])$. 
There is one singular point of convergence, at which a Burns 
metric bubbles off.  
\end{theorem}
\begin{remark}{\em
It is possible that $r$ could be taken to be infinite 
in the above theorem, but this would 
require a much more involved estimate of the Yamabe invariant.}
\end{remark}
\subsection{Acknowledgements} 
The author would like to thank first and foremost, Claude 
LeBrun, for originally suggesting this problem and 
providing extremely helpful comments along the way. 
Peter Kronheimer and Clifford Taubes also provided 
very helpful suggestions early on. Nobuhiro Honda provided 
crucial assistance in understanding the conformal 
geometry of LeBrun metrics. 
The author held several valuable discussions with Kazuo Akutagawa regarding 
the orbifold Yamabe problem.  
Denis Auroux, Simon Donaldson, Yat-Sun Poon, and Gang Tian 
also gave the author some important remarks that proved to be very 
useful when completing this work. 
Finally, thanks are given to the anonymous referee whose 
numerous remarks and suggestions greatly improved the 
exposition of the paper. 

\section{Monopole metrics}
\label{S1}
We first recall some basic definitions. 
\begin{definition}{\em
 A {\em{Riemannian 
orbifold}} $(M^n,g)$ is a topological space which is a 
smooth manifold of dimension $n$ with a smooth Riemannian metric 
away from finitely many singular points.  
At a singular point $p$, $M$ is locally diffeomorphic 
to a cone $\mathcal{C}$ on 
$S^{n-1} / G$, where $G \subset SO(n)$ 
is a finite subgroup acting freely on $S^{n-1}$.
Furthermore, at such a singular point, the metric is locally the 
quotient of a smooth $G$-invariant metric on $B^n$ under 
the orbifold group $G$.

 A {\em{Riemannian multi-fold}} $M$ 
is a connected space obtained from a finite collection of 
Riemannian orbifolds by finitely many identifications of points.
If there is only one cone at a singular point 
$p$, then $M$ is called {\em{irreducible}} at $p$, otherwise 
$M$  is called {\em{reducible}} at $p$. }
\end{definition}

We note that the notions of smooth orbifold,
orbifold diffeomorphism, and orbifold Riemannian 
metric are well-defined, see \cite{TV2} for 
background and references. Note also that our 
definition is very restrictive since we only 
allow isolated singular points. 

\begin{definition}{\em
 A smooth Riemannian manifold $(X^n,g)$ 
is called an asymptotically locally 
Euclidean (ALE) end of order $\tau$ 
if there exists a finite subgroup $G \subset SO(n)$ 
acting freely on $\mathbf{R}^n \setminus B(0,R)$ and a 
$C^{\infty}$ diffeomorphism 
$\Psi : X \rightarrow ( \mathbf{R}^n \setminus B(0,R)) / G$ 
such that under this identification, 
\begin{align}
g_{ij} &= \delta_{ij} + O( r^{-\tau}), \ \ \ 
\partial^{|k|} g_{ij} = O(r^{-\tau - k }),
\end{align}
for any partial derivative of order $k$ as
$r \rightarrow \infty$.  
A complete, noncompact Riemannian orbifold $(X,g)$ is called ALE 
if $X$ can be written as the disjoint union of a compact 
set and finitely many ALE ends. 
If all of the groups $G_j$ corresponding to 
the ends are trivial, then $(X,g)$ is called 
asymptotically flat (AF). } 
\end{definition}
For an integer $m \geq 1$, 
we let $\ZZ_m = \ZZ / m \ZZ \subset SU(2)$ be the cyclic 
group of matrices
\begin{align}
\label{su2}
\left(
\begin{matrix}
\exp^{2 \pi i p / m}   &  0 \\
0    & \exp^{-2 \pi i p / m } \\
\end{matrix}
\right),  \ \ 0 \leq p < m. 
\end{align}
acting on $\RR^4$, which is identified with $\CC^2$
via the map  
\begin{align}
(x_1, y_1, x_2, y_2) \mapsto (x_1 + i y_1, x_2 + i y_2) = (z_1, z_2).
\end{align}

\subsection{Gibbons-Hawking ansatz}
\label{GHansatz}
We briefly review the construction of 
Gibbons-Hawking multi-Eguchi-Hanson metrics,
from \cite{GibbonsHawking, Hitchin2}, as presented in 
\cite{AKL}. We also present a generalization to allow
orbifold points, by taking Green's functions with
integral weights.

Consider $\RR^3$ with the flat metric $g_{\RR^3} = dx^2 + dy^2 + dz^2$. 
Choose $n$ distinct points  $P = \{p_1, \dots, p_n \} \subset \RR^3$.
For each point $p_i$, we assign 
a multiplicity $m_i$, which is an 
integer satisfying $m_i \geq 1$, and 
let $N = \sum_{i = 1}^n m_i$ be the total multiplicity. 

Let $\Gamma_{p_j}$ denote the fundamental solution for the Euclidean
Laplacian based at $p_j$ with normalization 
$\Delta \Gamma_{p_j} = -2 \pi \delta_{p_j}$,
and let 
\begin{align}
V = \frac{1}{2}\sum_{i = 1}^n  m_i \Gamma_{p_i}.
\end{align}
Then $* dV$ is a closed $2$-form on $\RR^3 \setminus P$,
and $(1/ 2 \pi)[* dV]$ is an integral class in $
H^2 ( \RR^3 \setminus P, \ZZ )$.
Let $\pi: X_0 \rightarrow \RR^3 \setminus P$ 
be the unique principal $\U(1)$-bundle determined by the 
the above integral class.
By Chern-Weil theory, there is a connection form $\w \in H^1(X_0, \i \RR)$
with curvature form $\i (* dV)$. The Gibbons-Hawking metric is defined by 
\begin{align}
\label{GHmetric}
g_{\GH} =  V \cdot g_{\RR^3} - V^{-1} \w \odot \w.
\end{align}
Note the minus sign appears, since by convention our connection 
form is imaginary valued.
We define a larger manifold $X$ by attaching points $\tilde{p_j}$ 
over each $p_j$.
\begin{remark}{\em Choosing a different connection form will 
result in the same metric, up to diffeomorphism.}
\end{remark}
We summarize the main properties of $(X, g_{\GH})$ in 
the following proposition.  
\begin{proposition}
The Gibbons-Hawking multi-Eguchi-Hanson metric $(X_0, g_{\GH})$ 
extends to $X$ as a smooth Riemannian 
orbifold metric. 
At a point $\tilde{p}_i$ with multiplicity 
$m_i \geq 1$, $X$ has an orbifold structure 
with group $\ZZ / m_i \ZZ$, 
acting as in (\ref{su2}). 
The space $(X, g_{GH})$ is ALE with a single end of order $4$.
The group at infinity is the cyclic 
group $\ZZ / N \ZZ$ acting as in (\ref{su2}), 
where $N$ is the total multiplicity. 
\end{proposition}
\begin{proof}
The smooth case is 
discussed in \cite{AKL}, and is straightforward to 
adapt to the orbifold case. 
\end{proof}
This metric will be denoted by $g_{GH}(m_1 \cdot p_1, \dots, m_n \cdot p_n)$.  
Note that a small sphere around 
a $\ZZ_m$-orbifold point is diffeomorphic 
to the Lens space $L(m, 1)$, and
that $X$ is equipped with an isometric $S^1$ action, 
with fixed point set the finite set $\{\tilde{p}_1, \dots, \tilde{p}_n \}$.
\subsection{LeBrun hyperbolic ansatz}
\label{Lansatz}
We briefly review LeBrun's construction of K\"ahler
scalar-flat metrics on the blow-up of
$\CC^2$ at $n$ points on a line from \cite{LeBrunJDG}.
As in the Gibbons-Hawking case, we present 
a generalization which allows orbifold 
points. 

The LeBrun construction \cite{LeBrunJDG} is similar 
to the Gibbons-Hawking construction above, by replacing $\RR^3$ with the 
upper half-space model of hyperbolic space
\begin{align}
\H^3 = \{  (x,y,z) \in \RR^3,  z > 0 \},
\end{align}
with the hyperbolic metric $g_{\H^3} = z^{-2} ( dx^2 + dy^2 + dz^2)$. 
Choose $n$ distinct points $P = \{p_1, \dots, p_n \} \subset \H^3$.
For each point $p_i$, we assign 
a multiplicity $m_i$, which is an 
integer satisfying $m_i \geq 1$, and let 
$N = \sum_{i = 1}^n m_i$ be the total multiplicity.

Let $\Gamma_{p_j}$ denote the fundamental solution for the hyperbolic 
Laplacian based at $p_j$ with normalization 
$\Delta \Gamma_{p_j} = -2 \pi \delta_{p_j}$,
and let
\begin{align}
\label{Vdef}
V = 1 + \sum_{i = 1}^n  m_i \Gamma_{p_i}.
\end{align}
Then $* dV$ is a closed $2$-form on $\H^3 \setminus P$,
and $(1/ 2 \pi)[* dV]$ is an integral class in $
H^2 ( \H^3 \setminus P, \ZZ )$.
Let $\pi: X_0 \rightarrow \H^3 \setminus P$ 
be the unique principal $\U(1)$-bundle determined by the 
the above integral class.
By Chern-Weil theory, there is a connection form $\w \in H^1(X_0, \i \RR)$
with curvature form $\i (* dV)$. LeBrun's metric is defined by 
\begin{align}
\label{LBmetric}
g_{\LB} = z^2 (  V \cdot g_{\H^3} - V^{-1} \w \odot \w).
\end{align}
We define a larger manifold $X$ by attaching points $\tilde{p_j}$ 
over each $p_j$, and by attaching an $\RR^2$ at $z = 0$. 
The space $X$ is non-compact, and has the topology 
of an asymptotically flat space. Adding the point at 
infinity will result in a compact manifold $\hat{X}$. 
We summarize the main properties of $(X, g_{\LB})$ in 
the following proposition.  
\begin{proposition}[LeBrun \cite{LeBrunJDG}]
The metric $g_{\LB}$ extends to $X$ as a smooth
orbifold Riemannian metric. At a point $\tilde{p}_i$ with multiplicity 
$m_i \geq 1$, $X$ has an orbifold structure 
with group $\ZZ / m_i \ZZ$, 
acting as in (\ref{su2}). The space $(X, g_{\LB})$ is asymptotically flat 
K\"ahler scalar-flat with a single end of order $2$. 
By adding one point, this metric conformally compactifies to a
smooth self-dual 
conformal class on the compactification $(\hat{X}, [ \hat{g}_{\LB}])$.
If all points have multiplicity $1$, then $\hat{X}$ is diffeomorphic 
to $n \# \CP^2$.
\end{proposition}
\begin{proof} The smooth case is proved in \cite{LeBrunJDG}.
Furthermore, the case of one point taken with multiplicity $n$ 
was also considered in \cite[Section 5]{LeBrunJDG}, 
and the generalization to several points with 
multiplicity is straightforward. 
\end{proof}
The non-compact AF metric will be 
denoted by $g_{\LB}(m_1 \cdot p_1, \dots, m_n \cdot p_n)$,
while a metric on the compactification will 
typically be denoted by $\hat{g}_{\LB}$.
Note that a small sphere around a $\ZZ_m$-orbifold point is diffeomorphic 
to the Lens space $L(m, 1)$, and
that $\hat{X}$ is equipped with a conformal $S^1$ action, 
with fixed point set the finite set $\{\tilde{p}_1, \dots, \tilde{p}_n \}$
together with an $S^2$ corresponding to the boundary of $\mathcal{H}^3$. 
For $n = 0$, this construction gives the Euclidean metric 
on $\RR^4$, and for $n=1$, it yields the Burns metric, which 
conformally compactifies to the Fubini-Study metric on $\CP^2$. 
In other words, the Burns metric is the Green's function 
asymptotically flat space associated to the Fubini-Study metric. 

\subsection{Negative mass metrics}
\label{negativemass}
In \cite{LeBrunnegative}, LeBrun presented the first known examples of scalar-flat ALE 
spaces of negative mass, which gave counterexamples to extending the
positive mass theorem to ALE spaces. We briefly describe these as follows. Define
\begin{align}
g_{O\LB} = \frac{ dr^2}{ 1 + Ar^{-2} + B r^{-4}} +r^2 \Big[ \sigma_1^2 + \sigma_2^2
+ (  1 + Ar^{-2} + B r^{-4}) \sigma_3^2 \Big],
\end{align}
where $r$ is a radial coordinate, and $\{ \sigma_1, \sigma_2, \sigma_3 \}$ is a 
left-invariant coframe on $S^3 = {\rm{SU}}(2)$, and $A = n -2$, $B = 1 - n$. 
Redefine the radial coordinate to be $\hat{r}^2 = r^2 - 1$, 
and attach a $\CP^1$ at $\hat{r} = 0$. After taking a quotient by $\ZZ_{n}$, 
the metric then extends smoothly over this $\CP^1$, is ALE at infinity,
and is diffeomorphic to $\mathcal{O}(-n)$.  
The mass is computed to be $-4 \pi^2 (n -2)$, which is 
negative when $n > 2$. For $n =1$, this construction yields the 
Burns metric. For $n = 2$, this space is Ricci-flat, and 
is exactly the metric of Eguchi-Hanson. There is a close connection 
with the hyperbolic monopole metrics: the conformal compactification 
of these ALE spaces are confomal to  $\hat{g}_{LB}(n \cdot p_1)$,
a compactified LeBrun hyperbolic 
monopole orbifold metric with a single monopole point of multiplicity $n$
\cite[Section 5]{LeBrunJDG}.

\section{The Yamabe invariant}
\label{YI}

The Yamabe Problem asks whether there exists a 
conformal metric with constant scalar curvature on any closed Riemannian manifold,
and has been completely solved in the affirmative. 
We do not attempt to give a history of the Yamabe problem here, 
for this we refer the reader to \cite{Aubin, Schoen, LeeandParker}.
In what follows, let $(M,g)$ be a Riemannian manifold, and 
let $R$ denote the scalar curvature of $g$. 
Writing a conformal metric as $\tilde{g} = v^{\frac{4}{n-2}} g$, 
the Yamabe equation takes the form
\begin{align}
\label{Yam}
- 4 \frac{n-1}{n-2} \Delta v + R \cdot v = \lambda \cdot v^{\frac{n+2}{n-2}},
\end{align}
where $\lambda$ is a constant (note: we use the analyst's Laplacian). 
These are the Euler-Lagrange equations of the {\em Yamabe functional}, 
\begin{align}
\mathcal{Y}(\tilde{g}) =  Vol(\tilde{g})^{\frac{2-n}{n}} 
\int_M R_{\tilde{g}} dvol_{\tilde{g}},
\end{align}
for $\tilde{g} \in [g]$, where $[g]$ denotes the conformal class of $g$. 
An important related conformal invariant is the 
{\em Yamabe invariant} of the conformal class $[g]$: 
\begin{align}
Y(M,[g]) \equiv \underset{ \tilde{g} \in [g] }{\mbox{inf }} \mathcal{Y}(\tilde{g}).
\end{align}
In dimension $4$, Aubin's inequality states
\begin{align}
Y(M,[g]) \leq Y(S^4, [g_{S}]) = 8 \pi \sqrt{6},
\end{align}
with equality if $(M,g)$ is conformally 
equivalent to $(S^4, g_{S})$. 

\subsection{Upper estimate}
Next is an estimate from above on the 
Yamabe invariant of compactified LeBrun conformal classes. 
We begin with a short calculation.
\begin{lemma} 
\label{rholemma}
The function $\rho(p) = d_{\H^3} \big( p, (0,0,1) \big)$ 
satisfies
\begin{align} 
- 2 + \Delta \rho &= 4 e^{-2 \rho} V_1,
\end{align}
where
\begin{align}
V_1 &= 1 + \Gamma_{(0,0,1)} = 1 + \frac{1}{e^{2 \rho} -1}= \frac{1}{1 - e^{-2 \rho}}.
\end{align}
\end{lemma}
\begin{proof}
Since $\rho$ is the distance function 
of the hyperbolic metric, 
\begin{align}
\Delta \rho = 2 \coth \rho,
\end{align}
which yields
\begin{align} 
-2 + \Delta \rho &= 2 ( \coth \rho -1) = 
 2 \left( \frac{ e^{\rho} + e^{-\rho}}{ e^{\rho} - e^{- \rho}} - 1 \right)
= \frac{4}{ e^{2 \rho} - 1} = 4 e^{-2 \rho} V_1.
\end{align}
\end{proof}
For purposes of the Yamabe invariant, 
if $g$ is a $n$-pointed LeBrun metric, then 
we may assume that $p_1 = (0,0,1)$. To see this, 
apply a hyperbolic isometry $\phi$ to arrange so that 
$p_1 = (0,0,1)$. By \cite[Section 2]{HV}, 
there exists a lift of $\phi$ to $\Phi : X_0 \rightarrow X_0$ preserving 
the connection form. The metric $\Phi^* g$  will be conformal to 
the original, so will have the same Yamabe invariant.  

For $n=1$, the compactified LeBrun metric is conformal 
to $\CP^2$ with the Fubini-Study metric $g_{FS}$ \cite[Section 3]{LeBrunNN}.  
Consequently, 
\begin{align}
\label{YCP2}
Y(\hat{M},[\hat{g}_{LB}(p_1)]) = 12 \pi \sqrt{2}  = Y( \CP^2, [g_{FS}]). 
\end{align}
For $n > 1$ we have the following estimate of the Yamabe invariant. 
\begin{theorem}
\label{est1} 
Let $(\hat{M}, [\hat{g}_{LB}])$ be a Lebrun conformal class on $n \# \CP^2$
with $n > 1$. Without loss of generality assume that 
$p_1 = (0,0,1)$. Then  
\begin{align}
\label{yup2}
Y(\hat{M},[\hat{g}_{LB}]) <  12 \pi \sqrt{2} - \beta(p_2, \dots p_n), 
\end{align}
where  $\beta: \H^3 \times \cdots \H^3 \rightarrow 
\RR_+$, is a positive function, and 
approaches zero only if every point $p_2, \dots, p_n$ approaches
the boundary of $\H^3$. 
\end{theorem}
\begin{proof} 
We consider conformal changes of the form 
\begin{align}
\label{confform}
g = e^{2f} (  V \cdot g_H + V^{-1} \omega \odot \omega),
\end{align}
where $f : \H^3 \rightarrow \RR$. 
It is computed in \cite{LeBrunNN} that
\begin{align}
\label{LNNR}
R_g = 6 e^{-2f} V^{-1} ( -1 - \Delta f - |\nabla f|^2).
\end{align}
We will now take $f = - \rho$ as a test function in the Yamabe 
functional. The resulting metric $g$ will then be a smooth 
orbifold metric on the conformal compactification. 
Let $V = V_n$ correspond to an $n$-pointed LeBrun metric, with $p_1 = (0,0,1)$.
From Lemma~\ref{rholemma},
\begin{align}
R_g = 6 e^{2 \rho} V_n^{-1} ( \Delta \rho - 2) = 24 V_n^{-1} V_1
\end{align}
(since $|\nabla \rho| = 1$).
The Yamabe functional evaluated at $g$ is then 
\begin{align}
\mathcal{Y}(g) &= \int_M R_g d V_g \cdot 
\left( \int_M dV_g \right)^{-1/2}
= \int_M 24 V_n^{-1} V_1 dV_g\cdot 
\left( \int_M dV_g \right)^{-1/2}.
\end{align}
We next take a coordinate system $(x,y,z,\theta)$ where $\theta$
is an angular coordinate on the fiber for some trivialization. 
The volume element is
\begin{align}
 \sqrt{ \mbox{det}(g)} = e^{-4 \rho}(V^3 z^{-6} V^{-1})^{1/2} = e^{-4 \rho} V_n z^{-3}. 
\end{align}
In coordinates we then have
\begin{align}
\mathcal{Y}( g) = \int_M 24 V_1 e^{-4 \rho} z^{-3} dV_0 \cdot 
\left( \int_M V_n e^{-4 \rho} z^{-3}  dV_0 \right)^{-1/2},
\end{align}
where $dV_0 =  dx \wedge dy \wedge dz \wedge d \theta$.
Since the integrand is independent of $\theta$, 
we may integrate with respect to $\theta$ to obtain
\begin{align}
\label{rhofac}
\mathcal{Y}( g) = 24 \sqrt{2 \pi} \int_{\H^3} V_1 e^{-4 \rho} dV_{\H^3} \cdot 
\left( \int_{\H^3} V_n e^{-4 \rho} dV_{\H^3} \right)^{-1/2}.
\end{align}
Since $n \geq 1$, we must have $V_n \geq V_1$ (with strict 
inequality for $n >1$) and we obtain the estimate 
\begin{align}
\label{v1vn}
\mathcal{Y}( g) \leq 24 \sqrt{2 \pi} \left(  \int_{\H^3} V_1 e^{-4 \rho} dV_{\H^3} 
\right)^{1/2}. 
\end{align}
Using radial coordinates on $\H^3$, we compute that 
\begin{align}
\begin{split}
\label{v1int}
 \int_{\H^3} V_1 e^{-4 \rho} dV_{\H^3} 
&=  4 \pi \int_0^{\infty}  
\frac{1}{1 - e^{-2 \rho}} e^{-4 \rho} ( \sinh \rho)^2 d \rho\\
&= \pi  \int_0^{\infty} \left( e^{-2 \rho} - e^{-4 \rho} \right)  
d \rho = \frac{\pi}{4}.
\end{split}
\end{align}
Substituting this into \eqref{v1vn}, we obtain
\begin{align}
\mathcal{Y}( g) \leq 12 \pi \sqrt{2}= Y(\CP^2, [g_{FS}[). 
\end{align}
 If $n > 1$, this inequality is strict, and 
the only way it can be close to 
saturation is if all the points $p_i$ are
close to the hyperbolic boundary, since the only 
inequality used was $V_1 < V_n$. The existence 
of the function $\beta$ follows easily. 
\end{proof}

\section{The orbifold Yamabe invariant}
\label{OYP}

The orbifold Yamabe invariant of an orbifold conformal class is defined as 
in the smooth case: 
\begin{align}
Y_{orb}(M,[g]) = \inf_{\tilde{g} \in [g]} Vol( \tilde{g})^{\frac{2-n}{n}}
\int_M R_{\tilde{g}} dV_{\tilde{g}}. 
\end{align}
The analogue of Aubin's estimate 
and basic existence result is as follows.
\begin{theorem}[Akutagawa-Botvinik \cite{AB2}, Akutagawa \cite{Akutagawacoyi}] 
\label{t1}
Let $(M,g)$ be a Riemannian orbifold with 
singular points $\{p_1, \dots, p_k \}$, 
with orbifold groups $G_i \subset{\rm{SO}}(n)$, 
$i = 1 \dots k$. Then 
\begin{align} 
\label{orbest}
Y_{orb}(M,[g]) \leq  Y(S^n) \min_i |G_i|^{-\frac{2}{n}}. 
\end{align} 
Furthermore, if this inequality is strict, 
then there exists a smooth conformal metric $\tilde{g} = u^{\frac{4}{n-2}} g$
which minimizes the Yamabe functional (and thus has constant 
scalar curvature). 
\end{theorem}
We note that for the football $M = S^4/ \ZZ_n$ (as defined in 
the introduction) with the round metric $g_S$, we have 
\begin{align}
 Y_{orb}( M, [g_S])=  
Y(S^4, [g_S]) |n|^{-\frac{1}{2}} = \frac{ 8 \pi \sqrt{6}}{\sqrt{n}}.
\end{align}

 Given an orbifold with non-negative scalar curvature, 
one can use the Green's function for the conformal 
Laplacian to naturally 
associate with any point a scalar-flat ALE orbifold by 
\begin{align}
(M \setminus \{p\}, g_p = \Gamma_p^{ \frac{4}{n-2}}g ).  
\end{align}
An ALE coordinate system arises from using inverted 
normal coordinates in the metric $g$ in a neighborhood
of the point $p$. If we choose $p$ to be $p_k$, one of the 
orbifold points, then the end of this ALE space will 
correspond exactly to the group $G_k$. 

 The positive mass theorem does not hold in general 
 for ALE spaces as illustrated by LeBrun's negative mass examples 
discussed in Section \ref{negativemass}. Nakajima has proved a version 
of the positive mass theorem for spin ALE spaces 
with group $G \subset {\rm{SU}}(2)$ \cite{Nakajima},
in which the zero mass spaces are exactly the hyperk\"ahler 
ALE spaces classified by Kronheimer \cite{Kronheimer}.
This makes the orbifold Yamabe problem more
subtle than in the smooth case. Indeed, Schoen's test 
function from \cite{Schoen} will not prove strict inequality 
in \eqref{orbest} if the mass is non-positive. 

By an {\em{orbifold compactification}} of an ALE space $(X,g)$,
we mean choosing a conformal factor $u : X \rightarrow \RR_+$ 
such that $u = O(r^{-2})$ as $r \rightarrow \infty$. 
The space $(X, u^2 g)$ then compactifies to a
$C^{1,\alpha}$ orbifold. The next result states that the ALE spaces 
we will consider have {\em{smooth}} orbifold compactifications
with strictly positive orbifold Yamabe invariant.
\begin{proposition}
\label{posprop}
Let $(X,g)$ be either (i) a LeBrun hyperbolic monopole 
orbifold AF metric, or (ii) a Gibbons-Hawking 
orbifold ALE metric. Then there exists a 
$C^{\infty}$-orbifold conformal compactification 
$(\hat{X}, \hat{g})$ which satisfies $Y_{orb}([\hat{g}]) > 0$. 
\end{proposition}
\begin{proof} 
The existence of a smooth orbifold 
compactification follows directly 
from \cite[Proposition 12]{ChenLeBrunWeber},
the proof of which is based on twistor theory. 
A second proof, not using twistor theory, is obtained 
by locally solving the negative Yamabe problem 
near the orbifold point (which is a convex variational 
problem; this is solvable in the orbifold setting), and applying the 
removable singularity theorem for constant scalar curvature self-dual metrics 
\cite[Theorem 6.4]{TV2}. 

Next, one may find a conformal metric on the 
orbifold compactification whose scalar curvature does not 
change sign \cite[Lemma 3.4]{AB2}.
The strict positivity of the scalar curvature then 
follows using the strong maximum 
principle, as in \cite[Proposition 13]{ChenLeBrunWeber}.
Thus we must have  $Y_{orb}([\hat{g}]) > 0$.
We remark that it is not difficult to write down an explicit conformal 
factor on the compactification which has positive scalar curvature, 
but we leave this as an exercise for the interested reader.
\end{proof}
We next have an estimate for the Yamabe invariant 
of LeBrun orbifold metrics. 
\begin{theorem} 
\label{oyest}
Let $(\hat{M}, \hat{g}_{\LB}( m_1 \cdot p_1, \dots, m_n \cdot p_n) \big)$
be a conformally compactified LeBrun metric  
with total multiplicity $N = m_1 + \dots + m_n$. Without loss
of generality, assume that $p_1 = (0,0,1)$, 
and assume that all monopole points are contained in $B(p_1,r)$. Then 
\begin{align}
Y_{orb}([\hat{g}_{\LB}]) \leq \frac{12 \pi \sqrt{6}}{ \sqrt{N +2}} + O(r),
\end{align}
as $r \rightarrow 0$. 
\end{theorem}
\begin{proof}
As in the proof of Theorem \ref{est1}, we let $f = - \rho$, 
to obtain 
\begin{align}
\label{rhofac1}
\mathcal{Y}( g) = 24 \sqrt{2 \pi} \int_{\H^3} V_1 e^{-4 \rho} dV_{\H^3} \cdot 
\left( \int_{\H^3} V_n e^{-4 \rho} dV_{\H^3} \right)^{-1/2}.
\end{align}
Using \eqref{v1int}, we obtain 
\begin{align}
\label{rhofac2}
\mathcal{Y}( g) = 6 \sqrt{2} \pi^{3/2} \cdot 
\left( \int_{\H^3} V_n e^{-4 \rho} dV_{\H^3} \right)^{-1/2}.
\end{align}
Clearly, 
\begin{align}
V_N  = 1 +  \frac{N}{e^{2\rho} -1} + O(r),
\end{align}
as $r \rightarrow 0$. 
Using radial coordinates, we calculate
\begin{align*}
\int_{\H^3} \left( 
1 +  \frac{N}{e^{2\rho} -1} \right) e^{-4 \rho}  dV_{\H^3}
&= 4 \pi \int_0^{\infty}  \left( 
1 +  \frac{N}{e^{2\rho} -1} \right) e^{-4 \rho} (\sinh \rho )^2 d\rho \\
& = \frac{\pi}{12} \Big( 
-6 e^{-2t} - 3(N-2) e^{-4t} + 2(N-1) \Big) \Big|_0^{\infty} \\
&= \frac{(N + 2) \pi}{12},
\end{align*}
This yields
\begin{align}
\mathcal{Y}( g) &= 6 \sqrt{2} \pi^{3/2} \left( \frac{N + 2}{48} \right)^{-1/2} + O(r)
= \frac{12 \pi \sqrt{6}}{ \sqrt{N +2}} + O(r),
\end{align}
as $r \rightarrow 0$. 
\end{proof}
\begin{corollary} 
\label{orbex}
Let $(\hat{M}, \hat{g}_{LB})$ be as in 
Theorem \ref{oyest}, and assume that the highest multiplicity 
at any point is strictly less than $4(N+2)/9$.
Then there exists a radius $r > 0$ such that if 
all points $p_1, \dots, p_n \in B(p_1,r)$, then  
there exists 
a solution of the orbifold Yamabe problem on  $(\hat{M}, [\hat{g}_{LB}])$.
\end{corollary}
\begin{proof}
Let $m$ be the greatest integer strictly less than  $4(N+2)/9$,
and let $G$ be the cyclic group of order $m$, acting on 
$S^4 \subset \RR^5 = \RR^4 \times \RR^1$ as in \eqref{su2}. 
We have
\begin{align}
\frac{12 \pi \sqrt{6}}{ \sqrt{N+2}} =  \frac{8 \pi \sqrt{6}}{ \sqrt{ 4 (N+2)/9}}
<  \frac{8 \pi \sqrt{6}}{\sqrt{m}} = Y_{orb}(S^4/G, [g_{S}]).
\end{align}
Therefore, for $r$ sufficiently small, by Theorem \ref{oyest}, 
\begin{align}
\label{lthmest}
Y_{orb}(\hat{M},[\hat{g}_{LB}]) \lesssim \frac{12 \pi  \sqrt{6}}{\sqrt{N+2}}
< Y_{orb}(S^4/G, [g_{S}]).
\end{align}
If the highest multiplicity of any orbifold point is $m$, we see that 
that estimate \eqref{orbest} will be satisfied for 
$r$ sufficiently small, and 
Theorem \ref{t1} then yields a solution of the orbifold Yamabe problem. 
\end{proof}

\begin{remark}{\em
For $N = 3$, the highest multiplicity allowed is $[20/ 9] = 2$. 
This allows multiplicity $2$ points (but not a single multiplicity 
$3$ point). 
This proves the existence statement in Theorem \ref{imaint}.  
The convergence statements in Theorem 
\ref{imaint} will be proved later in Section \ref{convergence}. }
\end{remark}

\subsection{Proof of Theorem \ref{nonexist}}
First, adding a point $p$ at infinity, 
there exists a smooth orbifold conformal compactification 
$(\hat{X}, \hat{g})$ of $(X_n,g)$ by Proposition \ref{posprop}. 
Assume by contradiction that $\hat{g}$ is a constant scalar curvature 
metric on the compactification $\hat{X}$ in this conformal class.
Letting $E$ denote the traceless Ricci tensor, we recall the 
transformation formula:  if $g = \phi^{-2} \hat{g}$, then 
\begin{align}
E_g = E_{\hat{g}} + (m-2) \phi^{-1} \big( \nabla^2 \phi - (\Delta \phi/m) \hat{g} \big),
\end{align} 
where $m$ is the dimension, and the covariant derivatives are taken 
with respect to $\hat{g}$.
Since $g$ is Ricci-flat and $m=4$, we have 
\begin{align}
E_{\hat{g}} = \phi^{-1} \big( - 2\nabla^2 \phi + (\Delta \phi/2) g \big).
\end{align} 
We next use the argument of Obata \cite{Obata}; integrating on $\hat{X}$, 
\begin{align}
\begin{split}
\int_{\hat{X}} \phi | E_{\hat{g}}|^2 d\hat{V} & = \int_{\hat{X}} \phi E_{\hat{g}}^{ij} 
\left\{  \phi^{-1} \big( -2 \nabla^2 \phi + (\Delta \phi/2) g \big)_{ij} \right\}  d\hat{V}\\
& =   -2\int_{\hat{X}}  E_{\hat{g}}^{ij} \nabla^2 \phi_{ij} d\hat{V}
=-2 \lim_{\epsilon \rightarrow 0} \int_{\hat{X} \setminus B(p, \epsilon)}   
E_{\hat{g}}^{ij} \nabla^2 \phi_{ij} d\hat{V}.
\end{split}
\end{align}
Since $g$ is the Green's function metric associated to $\hat{g}$ at $p$, we have
\begin{align}
g = \phi^{-2} \hat{g} = G^{\frac{4}{m-2}} \tilde{g} \sim r^{-4} \tilde{g},
\end{align}
which implies that $\phi \sim r^{2}$ where $r$ is the distance to $p$ with respect 
to the metric $\hat{g}$. 
Continuing the above calculation, integration by parts yields 
\begin{align}
\int_{\hat{X}} \phi | E_{\hat{g}}|^2 d\hat{V} 
& = -2 \lim_{\epsilon \rightarrow 0} \left( \int_{ \partial  B(p, \epsilon)} 
E_{\hat{g}}^{ij} (\nabla \phi)_i \nu_j d \sigma
- \int_{\hat{X} \setminus B(p, \epsilon)}  (\nabla_j E_{\hat{g}}^{ij} \cdot \nabla_i \phi) d\hat{V} \right).
\end{align}
The second term on the right hand side is zero since the scalar curvature of $\hat{g}$ is constant 
(by the Bianchi identity), and the first term on the right hand side 
limits to zero since the integrand is bounded. Indeed, since $\hat{g}$ is a smooth 
orbifold, the curvature is bounded near $p$, and $|\nabla \phi| \sim r$ near $p$.
Consequently, $E_{\hat{g}} \equiv 0$, and $\hat{g}$ is Einstein. 

Since we have two Einstein metrics in the conformal class, the complete 
manifold $(X, g)$ admits a nonconstant solution of the equation
\begin{align}
\nabla^2 \phi = \frac{\Delta \phi}{m} g. 
\end{align} 
Such a solution is called a {\em{concircular scalar field}}, 
and complete manifolds which admit a non-zero solution were
classified by Tashiro \cite{Tashiro} (see also \cite{Kuhnel}), 
who showed that $(X,g)$ must be conformal to one of the following:
(A) a direct product $V \times J$, where $V$ is an 
$(m-1)$-dimensional complete Riemannian manifold
and $J$ is an interval, (B) hyperbolic space $\mathcal{H}^m$, 
or (C) the round sphere $S^m$.

 The hyperk\"ahler ALE spaces under consideration have 
second homology generated by embedded $2$-spheres
with self-intersection $-2$, with intersection matrix given by 
the corresponding Dynkin diagram \cite{Kronheimer}. 
If such a space were diffeomorphic to a product 
$V^3 \times J$, then any of the above spherical 
generators in $H_2$ would be homologous 
to a cycle in $V^3$, and would therefore have 
zero self intersection since such a cycle can 
be deformed to a disjoint cycle by translating it in 
the $J$ direction.
Cases (B) and (C) obviously cannot happen since 
$(X_n,g)$ is not locally conformally flat for $n > 1$. 
This is a contradiction, and the nonexistence is proved.
Finally, the non-existence of a solution, together with 
Theorem \ref{t1}, imply that the orbifold Yamabe invariant is 
maximal
\begin{align}
Y_{orb}(\hat{X}_n, [\hat{g}]) = \frac{ 8 \pi \sqrt{6}}{\sqrt{n}}. 
\end{align} 
This completes the proof of Theorem \ref{nonexist}.
\begin{remark}{\em
We point out that the Obata portion of the 
above proof does not hold if instead the compact manifold 
is assumed to be Einstein. For example, consider 
$\CP^2$ with the Fubini-Study metric $g_{FS}$, which is Einstein. 
The associated Green's function ALE space at any point is the Burns 
metric, which is scalar-flat but not Ricci-flat.}
\end{remark}
\begin{remark}{\em
It is clear from the above proof that 
Theorem \ref{nonexist} also holds for Ricci-flat ALE spaces in other 
dimensions, as long as they are not homeomorphic to a product $V^{m-1} \times J$,
and not locally conformally flat. }
\end{remark}
We conclude this section by noting that the proof of 
Theorem~\ref{nonexist} is also valid in case $X$ has non-trivial orbifold points. 
\begin{theorem}
\label{nonexist2}Theorem \ref{nonexist} holds if $(X_n,g)$ is a 
Gibbons-Hawking multi-Eguchi-Hanson orbifold. 
\end{theorem}
\begin{proof}
As shown in Proposition \ref{posprop}, there exists a smooth 
conformal compactification $(\hat{X},g)$. 
The proof of Theorem \ref{nonexist} above shows that 
any constant scalar curvature metric on $\hat{X}$ 
conformal to $\hat{g}$ must be Einstein. 
Consequently, there exists a concircular scalar field 
on $(X,g)$. An examination of Tashiro's proof shows 
that there are no orbifolds with isolated 
singularities in case (A), since the 
product with an interval would create at least a 
$1$-dimensional singular set. 
Cases (B) and (C) cannot occur either since these are 
locally conformally flat. 
\end{proof}

\section{Symmetric metrics}
\label{symsec}
We begin with some elementary hyperbolic geometry.
Fix the point $p_0 = (0,0,1) \in \mathcal{H}^3$,
and let $\rho_0(\cdot) = d(p_0, \cdot)$ denote the 
hyperbolic distance to $p_0$.  
Let $u : \H^3 \rightarrow \RR_+$ be defined by 
\begin{align}
\label{udef}
u =  \frac{\mbox{sech}(\rho_0)}{z} = \frac{2}{ 1 + x^2 + y^2 + z^2}. 
\end{align}
To see the second equality in \eqref{udef}, recall the following formula 
\begin{align}
\cosh \Big( d_{\mathcal{H}} ( p_1, p_2) \Big) = 1 + \frac{|p_1 - p_2|^2}{2 z_1 z_2}, 
\end{align}
where $p_i = (x_i, y_i, z_i)$, and the norm on the right is the
Euclidean norm (\cite[Theorem 4.6.1]{Ratcliffe}). 
From this, we obtain  
\begin{align}
\cosh \Big( \rho_0(p) \Big)
= 1 + \frac{ x^2 + y^2 + (z-1)^2}{2 z}
= \frac{ x^2 + y^2 + z^2 + 1}{2z}.
\end{align}  
\begin{lemma}
\label{ueqn}
The function $u$ satisfies the equation 
\begin{align}
\Delta_{Euc} u  + z^{-1}\partial_z u = - 2 u^3.
\end{align}
\end{lemma}
\begin{proof}
We define $\tilde{g} = u^2 ( dx^2 + dy^2 + dz^2 + z^2 d \theta^2 )$.
Letting $\hat{x} = z \cos \theta, \hat{y} = z \sin \theta$,
we obtain
\begin{align}
\tilde{g} = \frac{ 4 }{ (1 + x^2 + y^2 + \hat{x}^2 + \hat{y}^2)^2} 
( dx^2 + dy^2 + d \hat{x}^2 + d \hat{y}^2 ).
\end{align}
The right hand side is the spherical metric on $S^4$ in coordinates 
arising from stereographic projection. Consequently, 
from \eqref{Yam}, 
\begin{align}
\Delta u = -2 u^3. 
\end{align}
Writing out the Laplacian in the $(x,y,z,\theta)$-coordinates, we 
obtain 
\begin{align}
-2 u^3 = \Delta u = \frac{1}{z} \sum_i \partial_i ( z \partial_i u )
= \Delta_{Euc} u + z^{-1}\partial_z u.
\end{align}
\end{proof}

\subsection{Symmetries for ${\bf{n=2}}$}
\label{sn2}
In this subsection, we will only consider the LeBrun construction with 
two monopole points, $p_1$ and $p_2$.  Without loss of generality, by
applying a hyperbolic isometry, we may assume the points
$p_1 = (0,0,r_0),$ and $p_2 = (0, 0, r_0^{-1})$. 
It was shown in \cite{HV} that the automorphism 
group is 
\begin{align}
G \equiv ({\rm{U}}(1) \times {\rm{U}}(1)) \ltimes D_4, 
\end{align}
where $D_4$ is the dihedral group of order $8$. 
There is the index $2$ subgroup given by 
\begin{align}
K = ({\rm{U}}(1) \times {\rm{U}}(1)) \ltimes (\mathbb{Z}_2 \oplus \mathbb{Z}_2), 
\end{align}
which are exactly the lifts of hyperbolic isometries preserving 
the set of $2$ monopole points. 

For a metric to be $K$-invariant, it must in particular be 
invariant under the bundle $\U(1)$-action. So we consider 
only metrics of the form $v^2 g_{LB}$, where $v : H^3 \rightarrow \RR_+$.
We refer the reader to \cite[Section 2]{HV} for the details on 
lifting isometries of $\H^3$ to automorphisms of LeBrun metrics. 
\begin{proposition} 
\label{liso}
If the conformal automorphism $\Phi$ of $g_{LB}$ is the
lift of an isometry of hyperbolic space $\phi$, 
then it is an isometry of the metric $v^2 g_{LB}$ provided that 
 \begin{align}
 ( z \circ \phi)^2 \cdot (v \circ \phi)^2  = z^2 \cdot v^2.
\end{align}
\end{proposition}
\begin{proof} 
Using \eqref{LBmetric},
\begin{align}
{\Phi}^* ( v^2 g_{LB}) &= (v \circ \phi)^2 {\Phi}^* g_{LB}= 
(v \circ \phi)^2 \left(\frac{z \circ \phi}{z} \right)^2 g_{LB}
= v^2 g_{LB}.
\end{align}
\end{proof}
Next, define the metric
\begin{align}
\label{smetric}
\tilde{g}_{LB} & = u^2 \cdot g_{LB},
\end{align}
where $u$ is defined in \eqref{udef}. 
Let $\phi$ denote inversion in the unit sphere
\begin{align}
\phi(x, y, z) = \frac{1}{x^2 + y^2 + z^2} (x, y, z). 
\end{align}
\begin{proposition}
\label{phiiso}
The map ${\Phi}$ acts as an isometry of $\tilde{g}_{LB}$. 
\end{proposition}
\begin{proof}
We check 
\begin{align*}
( z \circ \phi)^2 (u \circ \phi)^2
 &= \frac{z^2}{\left(x^2 + y^2 + z^2\right)^2}  \cdot 
\frac{4}{\left( 1 + \frac{1}{x^2 + y^2 + z^2} \right)^2}
= \frac{4 z^2}{\left( 1 + x^2 + y^2 + z^2 \right)^2} = z^2 u^2,
\end{align*}
so the result follows from Proposition \ref{liso}.
\end{proof}
\begin{theorem}
\label{symesti}
For the LeBrun metric with $2$ monopole points 
$p_1 = (0,0, r_0), p_2 = (0,0, r_0^{-1})$, 
with $r_0 > 1$, the $K$-symmetric Yamabe invariant satisfies 
\begin{align}
Y_{K}(\hat{M},[g]) < 8 \sqrt{6} \pi - \beta(r_0),
\end{align}
where $\beta: (1, \infty) \rightarrow (0,\infty)$ 
satisfies  $\beta(r_0) \rightarrow 0$ as $r_0 \rightarrow \infty$. 
\end{theorem}
\begin{proof}
The identity component $\U(1) \times \U(1)$ is generated 
by (the lifts of) rotations around the $z$-axis
and the $\U(1)$ fiber rotation \cite[Proposition 2.14]{HV}.  
By Proposition \ref{liso}, these are isometries of the metric 
$\tilde{g}_{LB}$ defined in \eqref{smetric}. 
The group $K$ is generated by the identity 
component, by the lift ${\Phi}$ of
the inversion $\phi$, together with a lift of any 
reflection in the $(x,y)$-plane. The lift ${\Phi}$ acts as 
an isometry by Proposition \ref{phiiso}. 
The lift of a reflection in the $(x,y)$-plane
is also an isometry by Proposition \ref{liso}.
Consequently, $\tilde{g}_{LB}$ is a $K$-invariant 
metric. We next compute its Yamabe energy.

Take a coordinate system $(x,y,z,\theta)$ where $\theta$
is an angular coordinate system on the fiber for some 
trivialization. The volume element of $g_{LB}$ is 
\begin{align}
 \sqrt{ \mbox{det}(g)} = (V^3 z^2 V^{-1})^{1/2} = V z. 
\end{align}
Since $u$ depends only upon the $(x,y,z)$ coordinates, 
using Lemma \ref{ueqn}, 
we have for the Laplacian with respect to $g_{LB}$, 
\begin{align}
\Delta_{LB} u &= \frac{1}{Vz} \partial_i ( V^{-1} u_i V z) 
= V^{-1} (  \Delta_{Euc} u  + z^{-1}u_z) = V^{-1} ( - 2 u^3).
\end{align}
Since $g_{LB}$ is scalar-flat, from (\ref{Yam}) we have
\begin{align}
\tilde{R} = -6 u^{-3} \Delta u = 12 V^{-1}.
\end{align}
The Yamabe energy in coordinates is then given by
\begin{align}
\mathcal{Y}( \tilde{g}_{LB}) = \int_U z u^4 dV_0\cdot 
\left( \int_U z V u^4 d V_0 \right)^{-1/2},
\end{align}
where $d V_0 = dx \wedge dy \wedge dz \wedge d \theta$ is the 
coordinate volume element, and the region of integration 
is $U = \mathcal{H}^3 \times (0, 2 \pi)$. 
Since $n \geq 1$, we must have $V > 1$, and we obtain
the estimate 
\begin{align}
\mathcal{Y}( \tilde{g}_{LB}) < 
\left( \int z  u^4 d V_0 \right)^{+1/2} = 
8 \sqrt{6} \pi = Y(S^4, [g_{S}]).
\end{align}
The middle equality follows since for $n=0$, $\tilde{g}$ is the 
spherical metric, as seen in the proof of Lemma \ref{ueqn}. 
The existence of the function $\beta$ follows easily. 
\end{proof}

\subsection{Proof of Theorem \ref{nonsym}}
\label{nonex}
Recall the definition of LeBrun's negative mass metrics on $\mathcal{O}(-n)$
from Section \ref{negativemass} above. 
Since $g_{O\LB}$ is scalar-flat, the Yamabe equation for a metric 
$\hat{g} = f^{2} g$, $f > 0$, is 
\begin{align}
\label{Yam2}
- 6 \Delta f  = \lambda \cdot f^3,
\end{align}
where $\lambda > 0$ is a constant
(recall that from Proposition \ref{posprop}, there is 
a smooth conformal compactification $(\hat{X}, \hat{g})$ with 
strictly positive Yamabe invariant). 
We are interested in solutions which yield a smooth constant scalar curvature 
metric on the compactification.
In particular, we must 
have 
\begin{align}
\label{growth}
f = O(r^{-2}), \ |\nabla f| = O (r^{-3}),  |\nabla^2 f| = O(r^{-4}), 
\mbox{ as } r \rightarrow \infty. 
\end{align}
We are only interested in solutions for $r \in [1,\infty]$ which decay 
quadratically at $\infty$. Therefore, we make the change of 
coordinates $\hat{r}^2 = r^2 - 1$. In these new coordinates, the 
metric takes the form 
\begin{align}
g_{O\LB} = \left( \frac{ 1 + \hat{r}^2}{ n + \hat{r}^2} \right) d \hat{r}^2
+ (1 + \hat{r}^2) [ \sigma_1^2 + \sigma_2^2] + 
\frac{1}{1 + \hat{r}^2} \hat{r}^2 (n + \hat{r}^2) (\sigma_3^2 / n^2).
\end{align}
To obtain the actual K\"ahler scalar-flat metric on $\mathcal{O}(-n)$, 
one needs to attach a $\CP^1$ at
the origin, and quotient by $\ZZ_{n}$ (see \cite{LeBrunnegative}), but 
in the following we will just consider the metric to live on 
$\RR_+ \times S^3 \simeq \RR^4 \setminus \{0\}$, using 
$\hat{r}$ as radial coordinate.

As shown in \cite{LeBrunnegative}, the identity component of the 
isometry group of $g_{O\LB}$ is 
$\U(2)$. Obviously, 
any conformal factor for which the conformal metric 
is invariant under the subgroup
$\rm{SU}(2)$ must be radial. To yield a smooth metric on the 
compactification, $f$ must then satisfy the initial conditions
\begin{align}
f(0) = 1, \ f'(0) = 0. 
\end{align}

A computation shows that for radial $f$, the equation \eqref{Yam2} 
takes the form 
\begin{align}
\label{ryam}
\left( \frac{ n + (5 - 2n) \hat{r}^2}{\hat{r}(1 + \hat{r}^2)} \right) f_{\hat{r}}
+\left( \frac{ n + \hat{r}^2}{ 1 + \hat{r}^2} \right) f_{\hat{r} \hat{r}} 
= - \lambda f^3. 
\end{align}
If $n \geq 3$, then $5 - 2n < 0$. In this case, for $\hat{r}$ sufficiently 
large, the equation looks like 
\begin{align}
\label{signs}
(\text{negative})  f_{\hat{r}} + (\text{positive})  f_{\hat{r} \hat{r}} = (\text{negative}) f^3.
\end{align}
Recall that $ f \sim \hat{r}^{-2}$ for $\hat{r}$ large. Therefore, $f$ must be 
strictly decreasing for some $\hat{r}_0$ large, thus $f_{\hat{r}}(\hat{r}_0) < 0 $ . 
Examining the signs in \eqref{signs}, we see that $ f_{\hat{r} \hat{r}} (\hat{r}_0) < 0$. 
Consequently, the derivative of $f$ is strictly decreasing at $\hat{r}_0$. This implies that 
 $f_{\hat{r}}(\hat{r}) < 0 $ for all $\hat{r} \geq \hat{r}_0$, and therefore 
 $ f_{\hat{r} \hat{r}} (\hat{r}) < 0$  for all $\hat{r} \geq \hat{r}_0$.
This says that $f$ is concave, so $f$ must hit zero at some finite point, 
which is a contradiction. 

Finally, it is shown in \cite{LeBrunnegative}, that $g_{O\LB}$ for $n=2$ is isometric to 
the Eguchi-Hanson metric, which is hyperk\"ahler. Thus Theorem \ref{nonexist} 
can be applied to this case. This completes the proof of 
Theorem \ref{nonsym}.
\begin{remark}{\em
For $n =1$, it is easy to check that 
\begin{align}
f(\hat{r}) = \frac{1}{1 + \hat{r}^2} 
\end{align}
is a a solution of \eqref{ryam}. This is not a surprise, since for $n =1$, 
$g_{O\LB}$ is the Burns metric, 
and $\hat{g} = (1 + \hat{r}^2)^{-2} g_{O\LB} =  r^{-4} g_{O\LB}$ is the Fubini-Study metric \cite{LeBrunnegative}. }
\end{remark}

\section{Integral formulas}
\label{GBHS}
For an ALE space $X$ with several ends $E_i$, and orbifold 
singularity points $p_j$, we have the 
signature formula
\begin{align}
\label{signature}
\tau(X) = \frac{1}{12 \pi^2} \left( \int_X |W^+_g|^2 dV_g - \int_X
|W^-_g|^2 dV_g \right) - \sum_i \eta( S^3 / G_i )
+ \sum_j \eta( S^3 / G_j' ),
\end{align}
where $G_i \subset SO(4)$ is the group corresponding to
the $i$th end, $\eta(  S^3 / G_i)$ is the 
$\eta$-invariant, and $G_j'$ are the groups corresponding 
to the orbifold points $p_j$. The Gauss-Bonnet formula in this context is
\begin{align}
\label{GB}
\chi(X) = \frac{1}{8 \pi^2} \left(
\int_X |W_g|^2 dV_g 
+ 4\int_X \sigma_2 dV_g \right)
+ \sum_i \frac{1}{|G_i|}
+ \sum_j \left( 1 -  \frac{1}{|G_j'|} \right),
\end{align}
where 
\begin{align}
 4\int_X \sigma_2 dV_g = - \frac{1}{2}\int_X |E|^2 dV_g + \frac{1}{24} \int_X R^2 dV_g,
\end{align}
and $E$ is the traceless Ricci tensor. 
See \cite{Hitchin} for a nice discussion of these formulas.
In this section, we will compute these
for various examples. 
\subsection{Self-dual metric on $n \# \CP^2$}
In this case, we have $\chi(M) = n +2, \tau(M) = n$. Thus 
\begin{align}
12 \pi^2 n = \int_M |W^+_g|^2 dV_g,  
\end{align}
and
\begin{align}
8 \pi^2 (n + 2) &= \int_M |W_g^+|^2 dV_g + 4 \int_M \sigma_2 dV_g.
\end{align}
Combining these gives
\begin{align} 
\label{s2int}
4 \pi^2 ( 4 - n) = 4\int_M \sigma_2 dV_g. 
\end{align}
This yields an estimate 
for the Yamabe invariant of positive 
self-dual metrics on $2 \# \CP^2$ and $3 \# \CP^2$.
\begin{proposition}
\label{23prop}
Let $0 \leq n \leq 3$. 
If $g$ is a self-dual metric on $M = n \# \CP^2$ with 
positive scalar curvature, then 
\begin{align}
Y(M,[g]) \geq  4 \pi \sqrt 6 \sqrt{(4-n)}.
\end{align}
Thus for $n=2$, $Y(M,[g]) \geq 8 \pi \sqrt{3}$,
and for $n =3$, $Y(M,[g]) \geq 4 \pi \sqrt{6}$. 
\end{proposition}
\begin{proof}
Using the inequality $\sigma_2 \leq R^2/24$, we obtain
\begin{align}
4 \pi^2 (4-n) \leq  \frac{1}{24}  \int_M R^2,
\end{align}
Since the self-duality condition is conformally invariant, 
we may conformally change to a Yamabe minimizer and obtain 
\begin{align} 
4 \pi \sqrt 6 \sqrt{4-n} \leq Y(M, [g]). 
\end{align}
\end{proof} 
Restricting to the class of monopole metrics, we have the 
following. 
For $n = 0$, $\hat{g}_{\LB}$ is conformal to the round metric 
on $S^4$, so we have $Y(\hat{M},[\hat{g}]) = 8 \pi \sqrt{6}$.
For $n =1$,  $\hat{g}_{\LB}$ is conformal to the  Fubini-Study metric,
so we have we have $Y(\hat{M},[\hat{g}]) = 12 \pi \sqrt{2}$.
For $n =2$, we have the lower estimate on the 
Yamabe invariant stated in Theorem~\ref{n=2thm}.

We also make the following observation
\begin{proposition} 
If $g$ is a self-dual metric on $2 \# \CP^2$ 
or $3 \# \CP^2$ with 
positive scalar curvature, then $g$ is 
conformal to a metric with positive Ricci 
curvature. 
\end{proposition}
\begin{proof}
From \eqref{s2int}, we see that the 
conformal invariant $\int_M \sigma_2$
is positive when $n \leq 3$.  It follows from 
\cite{CGY1} that there exists a conformal metric with 
$R > 0$ and $\sigma_2 > 0$ pointwise 
(see also \cite{GurskyViaclovskyJDG}). 
Such a metric necessarily has positive 
Ricci curvature. 
\end{proof}
This was proved in \cite{LeBrunNN} for $n =2$, 
and for $n = 3$ under certain conditions on the 
$3$ monopole points. This was done by an explicit construction, 
as the result of Chang-Gursky-Yang was not known 
at the time.  An interesting fact is that positive 
Ricci metrics exist for {\em{any}} positive scalar curvature self-dual 
metric on $3 \# \CP^2$, not only those obtained by the 
LeBrun hyperbolic ansatz.

\subsection{Single monopole point with multiplicity}
We take a LeBrun metric with a single monopole point 
of multiplicity $n$, and compactify to a self-dual 
orbifold $\hat{M}$ with a single orbifold point of type $A_{n-1}$. 
It was shown in \cite{LeBrunJDG} that this is the same 
as the conformal compactification of LeBrun's ALE 
metrics on $\mathcal{O}(-n)$ found in \cite{LeBrunnegative}. 
The characteristic numbers are $\chi(\hat{M}) = 3, \tau(\hat{M}) = 1$,
and (see \cite{Nakajima}) 
\begin{align}
\eta( S^3 / G) = \frac{(n-1)(n-2)}{3n}. 
\end{align}
We thus have
\begin{align}
\label{signature2}
1 = \frac{1}{12 \pi^2} \int_{\hat{M}} |W^+_g|^2 dV_g  - \frac{(n-1)(n-2)}{3n},
\end{align}
(the minus sign is due to reversed orientation), so 
\begin{align}
 \int_{\hat{M}} |W^+_g|^2 dV_g = \frac{n^2 + 2}{3n} 12 \pi^2.
\end{align}
The Gauss-Bonnet formula yields 
\begin{align}
3 = \frac{n^2 + 2}{ 2n} 
+ \frac{1}{8 \pi^2} 4\int_{\hat{M}} \sigma_2 dV_g 
+ 1 -  \frac{1}{n},
\end{align}
which simplifies to 
\begin{align}
4 \pi^2 (4-n)  =  4\int_{\hat{M}} \sigma_2 dV_g.
\end{align}
\subsection{Gibbons-Hawking multi-Eguchi-Hanson}
For this example with $n$ points, we have $\tau(GH) = n-1 , \chi(GH) = n$.
Let $\hat{M}$ be the conformal compactification 
to a compact orbifold with a single $\ZZ_n$
singularity of type $A_{n-1}$. Reversing orientation, the metric 
is self-dual. We thus have $\tau(\hat{M}) = n - 1, \chi(\hat{M}) = n + 1$.
The signature formula is
\begin{align}
n-1 = \frac{1}{12 \pi^2} \int_{\hat{M}} |W^+_g|^2 dV_g 
+ \frac{(n-1)(n-2)}{3n},
\end{align}
which yields 
\begin{align}
\int_{\hat{M}} |W^+_g|^2 dV_g = \frac{(n+1)(n-1)}{n} 8 \pi^2. 
\end{align}
\begin{remark}{\em
Note that 
\begin{align}
\frac{(n+1)(n-1)}{n} 8 \pi^2 +  \frac{n^2 + 2}{3n} 12 \pi^2 = 12 n \pi^2,
\end{align}
which reflects the fact that the Lebrun metrics will degenerate
to a generalized connect sum of a LeBrun orbifold and 
a compactified multi-Eguchi-Hanson space as the 
$n$ monopole points tend to a single point, as observed in 
\cite{LeBrunJDG}. This will be made more precise in Section \ref{convergence}. 
}
\end{remark}
Returning to the example, 
the Gauss-Bonnet formula yields 
\begin{align}
\label{GB2}
n+1 =  \frac{(n+1)(n-1)}{n}
+ \frac{1}{8 \pi^2} 4\int_{\hat{M}} \sigma_2 dV_g 
+ 1 -  \frac{1}{n},
\end{align}
which simplifies to 
\begin{align}
\frac{2}{n}  8 \pi^2 = 4 \int_{\hat{M}} \sigma_2 dV_g.
\end{align}
This implies the inequality 
\begin{align}
\frac{2}{n}  8 \pi^2 \leq \frac{1}{24} \int_{\hat{M}} R^2 dV_g.
\end{align}
If there existed a Yamabe minimizing metric $\hat{g}$
in the conformal class, then 
\begin{align}
\frac{ 8 \pi \sqrt{6}}{\sqrt{n}} \leq Y_{orb}( \hat{M}, [\hat{g}]).
\end{align}
But the Akutagawa-Botvinnik inequality in Theorem \ref{t1} says the 
reverse, so we must have equality. Tracing through the inequalites, 
this says that $|E| = 0$.   This gives an alternative 
proof that any Yamabe minimizer would have to 
be Einstein  (and therefore cannot exist, recall the 
proof of Theorem \ref{nonexist} above). 

\section{Orbifold convergence}
\label{convergence}

\begin{figure}
\includegraphics{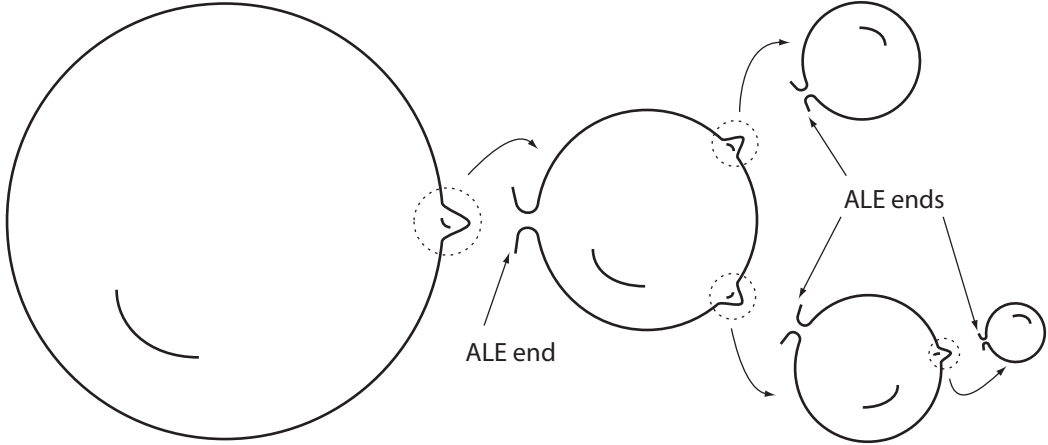}
\caption{A bubble tree with several levels of scaling. The dotted 
circles enclose the regions of singularity formation, and are
close to orbifold singularities (in the Cheeger-Gromov sense).
The curved arrows represent these regions viewed through a 
powerful microscope.}   
\label{bubbletree}
\end{figure}

We briefly describe the structure of ``bubble-trees''. 
Let $(M_i, h_i)$ be a sequence of metrics converging to an 
orbifold $(M_{\infty}, h_{\infty})$ in the Cheeger-Gromov
sense. 
At the first level of bubbling (the 
lowest level of rescaling), the 
ALE orbifolds $(X_{j}, g_{j})$, $j = 1 \dots k$,
bubble off, with each orbifold $X_j$ corresponding to singular 
points of convergence $p_j \in M_{\infty}, j = 1 \dots k$, of the 
original sequence $(M_i, h_i) \rightarrow (M_{\infty}, h_{\infty})$. 
Each $(X_j, g_{j})$, $j = 1 \dots k$, 
is a pointed rescaled Cheeger-Gromov limit of the original sequence,  
having singular points of convergence $p_{j, j'} \in X_{j}$, $j' = 1 \dots k_{j}$.
At each singular point $p_{j,j'}$, there is a further rescaling
which converges to an ALE orbifold $(X_{j, j'}, g_{j, j'})$, 
as above with singular points of convergence $p_{j, j', j''} \in X_{j,j'}$, 
$j'' = 1 \dots k_{j,j'}$. 
This procedure is then repeated and terminates in 
finitely many steps, see Figure \ref{bubbletree}. 
We refer the reader to \cite{TV2} for more details about
this procedure and further references.

\subsection{Bubble-tree structure for hyperbolic monopole metrics}
The following theorem describes the bubble formation
for compactified LeBrun metrics.
\begin{theorem}
\label{LBconv}
Let $(M,g_i)$ be a sequence of $n$-pointed LeBrun metrics with 
monopole points $\{ p_{i,1}, \dots, p_{i,n} \}$. 
Assume that as $i \rightarrow \infty$ that these
points converge to 
\begin{align}
\{ p_{i,1}, \dots, p_{i,n} \} \rightarrow
\{ m_1 \cdot p_{\infty,1}, \dots, m_k \cdot p_{\infty,k} \},
\end{align}
as $i \rightarrow \infty$ with $ p_{\infty,j} \in \H^3$
for $1 \leq j \leq k_1$ and $ p_{\infty,j} \in \partial \H^3$
for $k_1 < j \leq k$, allowing for multiplicity. 
Then there exist metrics $\hat{g}_j$ on the conformal 
compactification so that $(\hat{M}, \hat{g}_j)$ converges
to the compactified LeBrun orbifold
\begin{align}
(\hat{M}, \hat{g}_j) \rightarrow 
\big(M_{\infty}, \hat{g}_{\LB}(  
m_1 \cdot p_{\infty,1}, \dots, m_{k_1} \cdot p_{\infty, k_1})\big),
\end{align}
as $j \rightarrow \infty$ in the Cheeger-Gromov sense.  
There are finitely many bubbles, and the bubble-tree structure 
is as follows. For each subcollection 
of points limiting to a finite point with multiplicity greater 
than one, the bubble-tree structure is a tree
of Gibbons-Hawking multi-Eguchi-Hanson orbifold ALE spaces.  
For each subcollection of points limiting to a 
boundary point, then a LeBrun orbifold AF metric
is the first bubble at that point, with subsequent bubbles being 
Gibbons-Hawking orbifolds as in the previous case. 
The neck regions are modeled on annuli in Euclidean spaces 
$\RR^4 / \ZZ_k$, with the group action as in \eqref{su2}.
\end{theorem}
\begin{proof}
By a sequence of conformal transformations, we normalize 
the sequence so that $p_{i,1} = (0,0,1)$.
Choose a conformal factor $u : \H^3 \rightarrow \RR_+$ 
so that $\hat{g}_i = u^2 g_{i}$ is a sequence of 
smooth metrics on the compactification $\hat{M}$ as follows.
Viewing $\H^3$ as the upper half space, without loss of generality 
we may assume that the limiting boundary points lie in some 
compact set of the $(x,y,0)$-plane. We choose the compactifying 
conformal factor $u$ to be identically 
$1$ on some open neighborhood $U$ containing this compact set
and such that $U$ also contains all monopole points. 
Furthermore, choose $u$ to be asymptotic to \eqref{udef} on  
$\partial \H^3 \setminus U$, so that $\hat{g}_i$ is a  
smooth metric on the compactification.  
To understand the bubble structure, we can therefore ignore 
the compactifying conformal factor in the following argument.  

We first consider the case of several points limiting to a higher multiplicity 
point, say $m_1 \cdot p_{\infty,1}$ with $m_1 >1$. 
If all points limit to $p_{\infty,1}$ at 
a uniform rate, then we rescale the metric so that the 
points are minimally seperated by distance $1$. 
As in \cite[page 237]{LeBrunJDG}, this is equivalent 
to rescaling the hyperbolic metric to become the flat 
metric, and the function $V$ limiting to the sum of
Euclidean Green's functions (without a constant). Thus 
the rescaled limit is a Gibbons-Hawking metric. 
If the points do not limit to $m_1 \cdot p_{\infty,1}$ at a uniform 
rate, then we do the following.  Rescale the metric the smallest 
amount so that the maximum distance between these points is $1$. 
We will then see several ``clusters'' of points limiting to distinct 
points in a unit ball. Thus the limit will be a Gibbons-Hawking 
ALE orbifold. At each orbifold point $q$ of this limit, we return 
to the original sequence and rescale so that the maximum 
distance between the subcollection of points limiting to $q$ 
in the first rescaling is $1$. At this scaling, we will then see 
several new subclusters limiting to distinct points in a unit 
ball, so again we find Gibbons-Hawking orbifolds at a different scale. 

Next, for a subcollection of points limiting to a boundary 
point, by a conformal transformation, we arrange so that this cluster 
of points is contained in a unit ball around 
some finite point in $\H^3$, say $p= (0,0,1)$. The pointed limit
(based at $p$) of such rescaled metrics is a LeBrun hyperbolic monopole 
AF orbifold metric (since in this gauge, all other points will limit
to $\partial \mathcal{H}^3$. This is conformal related to the 
original, but it is easy to see that in the original scaling, 
precisely this AF orbifold metric bubbles off. 
We will illustrate this in a simple case, the general argument 
is the same. Consider the case of 2 monopole points. 
Let $p_{j,1} = (0,0,1),  p_{j,2} = (0,0,j^{-1})$. 
Let $\phi_j (x,y,z) = (jx, jy, jz)$, this is a hyperbolic
isometry. By \cite[Section 2]{HV}, there is a lift 
${\Phi}_j$ of $\phi$ preserving the connection 
form. We then have
\begin{align}
 {\Phi}_j^* ( g_{LB} ( p_{j,1}, p_{j,2}))
= j^2 \Big( g_{LB} \big( (0,0,j), (0,0,1) \big) \Big).
\end{align}
In other words, we see that a scaling of the original 
metric is isometric to another LeBrun metric.
Consequently, the bubble will be a LeBrun AF metric
(in this case a Burns metric). 
In general, there could be clusters of points tending 
towards $\partial \H^3$, and then clearly all deeper 
bubbles at such points will be Gibbons-Hawking 
orbifold metrics, as in the first paragraph. 

Since the number of monopole points is bounded by $n$, 
this procedure must terminate in finitely many 
steps, and thus there are finitely many {\em{non-trivial}} 
bubbles. A trivial bubble, or neck region, will 
arise at any intermediate scaling between 
a non-trivial ALE space, and the previous 
orbifold point onto which is it glued. All of the above 
ALE spaces are K\"ahler and have an ALE-coordinate system in which 
the group action is as in \eqref{su2}, so the structure 
of the neck regions is clearly as stated. 
\end{proof}
\begin{remark}{\em
To clarify, we explain the above procedure in 
two cases. Consider the case of $3$ points, and 
normalize so that $p_{j,1} = (0,0,1)$. Choose $p_{j,2}$ and $p_{j,3}$ 
to be at distance $j^{-2}$ from each other, but 
such that both of these are at distance $j^{-1}$ from $p_{j,1}$. 
In this case, the original limit is a LeBrun metric with a 
single multiplicity $3$ point. To see the first bubble, 
we rescale the picture so that the distance between 
$p_{j,1}$ and $p_{j,2}$ is $1$. Thus the first bubble will be a 
Gibbons-Hawking orbifold with one point of 
multiplicity $1$, and another point of multiplicity 
$2$. The second, deepest, bubble will rescale so that 
$p_{j,2}$ and $p_{j,3}$ are at distance $1$, and 
this bubble will be a Gibbons-Hawking metric 
with $2$ points of multiplicity $1$, which is none 
other than the classical Eguchi-Hanson metric. 

For the next example, consider again the case of 
$3$ points. Assume that $p_{j,1}  = (0,0,1)$, and choose 
$p_{j,2}$ and $p_{j,3}$ to limit to the boundary 
point $(0,0,0)$. The limit of the original sequence will 
be a LeBrun metric with a single monopole point
(which we know is conformal to the Fubini-Study 
metric).  To understand the bubbling, apply a conformal 
transformation so that $p_{j,2}$ and $p_{j,3}$ 
limit to $(0,0,1)$ (and then $p_{j,1}$ will limit to 
the boundary of hyperbolic space). The limit 
will now be a LeBrun AF metric with a single point 
of multiplicity $2$. This AF space will be the first 
bubble. Upon further rescaling to separate these two points, 
we find the second deepest bubble to 
be an Eguchi-Hanson metric. 
}
\end{remark}

\subsection{Sobolev constants and Yamabe invariants}
We proceed with the definition of the Sobolev constant. 
\begin{definition}{\em
 For $\hat{M}$ compact, we define the Sobolev constant $C_S$  
as the best constant $C_S$ so that 
for all $f \in C^{0,1}(\hat{M})$ we have
\begin{align}
\label{mainSob}
\Vert f \Vert_{L^{4}} \leq C_S \left(   
\Vert \nabla f \Vert_{L^2} + \Vert f \Vert_{L^2} \right).
\end{align}
For $M$ non-compact, $C_S$ is defined to be the best constant so that 
\begin{align}
\label{mainSob2}
\Vert f \Vert_{L^{4}} \leq C_S   
\Vert \nabla f \Vert_{L^2},
\end{align}
for all $f \in C^{0,1}(M)$ with compact support. 
}
\end{definition}
We next prove the lower estimate on the Yamabe invariant in Theorem \ref{maint}.
\begin{theorem} 
\label{maintp1}
Let $(\hat{M},\hat{g})$ be a LeBrun self-dual conformal class
on $n \# \CP^2$ with $n$ monopole points $\{ p_1, \dots, p_n\}$.
Assume that all monopole point are contained in a 
compact set $\mathcal{K} \subset \mathcal{H}^3$. Then
there exists a constant $\delta_n > 0$ depending 
only upon $n,\mathcal{K}$ such that
\begin{align}
\label{mest2}
0 < \delta_n \leq Y(\hat{M}, [\hat{g}]).
\end{align}
\end{theorem}
\begin{proof}
Let $\hat{g}_{\LB,i}$ be a sequence of such LeBrun metrics,
with monopole points $p_{i,j}$, $j = 1 \dots n$, 
and assume that
\begin{align}
\{ p_{i,1}, \dots, p_{i,n} \} \rightarrow
\{ m_1 \cdot p_{\infty, 1}, \dots, m_k \cdot p_{\infty, k} \},
\end{align}
as $i \rightarrow \infty$ with $ p_{\infty,j} \in \H^3$
for $1 \leq j \leq k$ allowing for multiplicity (from 
our assumption, there are no limit points on $\partial \H^3$).  
Inspired by \cite{LeBrunNN}, let 
\begin{align}
f_i = - ( \rho_{i,1} + \dots + \rho_{i,n})/n,
\end{align}
where $\rho_{i,j}(\cdot) = d_{\H^3}(p_{i,j}, \cdot)$.
Without loss of generality, by conformal invariance, 
we may assume that $p_{i,1} = (0,0,1)$. 
Consider the conformally compactified 
$\tilde{g}_i = e^{2f_i}z^{-2} g_{\LB,i}$, as in \eqref{confform}.  
Note that since all points are contained in a compact set
around $p_{i,1} = (0,0,1)$, $f_i$ is bounded from below 
on any compact set, and $f_i(x) \sim - p_{i,1}(x)$ as 
$x \rightarrow \partial \H^3$. 
This implies that $\lambda < Vol(\tilde{g}_i) < \Lambda$, 
for some positive constants $\lambda$ and $\Lambda$.
The same argument in Theorem \ref{LBconv}, 
then shows that $\tilde{g}_i$ converges to 
a non-trivial limiting orbifold $(\hat{M}_{\infty}, \hat{g}_{\infty})$ 
with finitely many ALE orbifold 
bubbles (our conformal factor now varies with $i$, but it remains 
bounded and strictly positive near the singularities, so 
the same argument applies). 
Since $|\nabla \rho_{i,j}| = 1$, this implies the inequality 
\begin{align}
|\nabla f_i|^2 <  1/n.
\end{align}
From \eqref{LNNR}, we estimate the scalar curvature
\begin{align*}
R_{\tilde{g}_i} &= 6 e^{-2f_i} V_i^{-1} ( -1 - \Delta f_i - |\nabla f_i|^2)\\
& > 6 e^{-2f_i} V_i^{-1} ( -2 - \Delta f_i)
=  6 e^{-2f_i} V_i^{-1} n^{-1} \sum_j (-2 + \Delta \rho_{i,j}),
\end{align*}
where $V_i$ is defined in \eqref{Vdef}, corresponding to the 
$n$ monopole points $p_{i,j}$.  
Using Lemma \ref{rholemma}, we obtain the estimate 
\begin{align}
R_{\tilde{g}_i} > 24 e^{-2f_i}  V_i^{-1} n^{-1} \sum_j e^{-2 \rho_{i,j}} V_{i,j},  
\end{align}
where $V_{i,j} = 1 + \Gamma_{p_{i,j}}$. For each $j$, the function 
$V_i^{-1} V_{i,j}$ is smooth and uniformly positive (with a lower 
bound independent of $i$), so we have 
\begin{align}
R_{\tilde{g}_i} > C e^{-2f_i} \sum_j e^{-2 \rho_{i,j}},
\end{align}
for some constant $C > 0$. 
Obviously, for some $j_0$ we must have 
\begin{align}
\frac{1}{n} \sum_j \rho_{i,j} \geq \rho_{i,j_0}.
\end{align}
Therefore, we have the estimate
\begin{align}
\label{Rdelta}
R_{\tilde{g}_i} > \delta > 0.  
\end{align}

 We also claim that there is a uniform bound on the Sobolev constant
\eqref{mainSob}. 
The limit space $(\hat{M}_{\infty}, \hat{g}_{\infty})$ 
has bounded Sobolev constant \eqref{mainSob} 
since it is a compact orbifold. Furthermore, all of the bubbles 
are ALE orbifolds, which have bounded Sobolev 
constants \eqref{mainSob2} by \cite{Bartnik}, and each are glued onto the previous 
orbifold by Euclidean neck regions. Therefore, using a standard partition 
of unity argument (at possibly several scales) and scale invariance 
of the Sobolev constant
\eqref{mainSob2}, it follows that the Sobolev constant of the sequence 
is uniformly bounded. This is proved in detail in 
\cite[Proposition 2.2]{JoyceCSC}, and the generalization to this case 
is straightforward, so we omit the details.  

To finish the proof, we include the following standard argument. 
Assume by contradiction that there is a sequence of unit-volume 
Yamabe minimizers in each conformal class $g_{Y,i} = v_i^{2} \tilde{g}_i$
with $R_i \rightarrow 0$ as $i \rightarrow \infty$, 
where $R_i$ is the (constant) scalar curvature of $g_{Y,i}$. 
The Yamabe equation is 
\begin{align}
-6 \Delta v_i + R(\tilde{g}_i) v_i = R_i v_i^3. 
\end{align}
Multiply by $v_i$ and integrate to get 
\begin{align}
6 \int |\nabla v_i|^2 dV_{\tilde{g}_i} 
+ \int R(\tilde{g}_i) v_i^2  dV_{\tilde{g}_i}  = \int R_i v_i^4  dV_{\tilde{g}_i}
= R_i Vol(g_{Y,i}) = R_i. 
\end{align}
The right hand side limits to zero,
therefore the left-hand side does also
as $i \rightarrow \infty$. From \eqref{Rdelta}, 
$R(\tilde{g}_i)$ is uniformly positive, so the $W^{1,2}$ norm of 
$v_i$ can be made arbitrarily
small. Using the uniform Sobolev inequality and lower volume bound 
for $\tilde{g}_i$, 
\begin{align}
1 = Vol( g_{Y,i})^{1/4} = \Vert v_i \Vert_{L^{4}} 
\leq C_S(\tilde{g}_i) \left( \Vert \nabla v_i \Vert_{L^2} 
+ \Vert v_i \Vert_{L^2} \right) \leq C' \Vert v_i \Vert_{W^{1,2}},
\end{align}
which is a contradiction for $i$ sufficiently large. 
\end{proof}
\subsection{Convergence of constant scalar curvature metrics}
We are now in a position to describe the possible limits 
of the Yamabe minimizers. The following Theorem implies the 
main part of Theorem \ref{maint}. 
\begin{theorem}
\label{poslim}
Fix $n$, and let $(M, g_i)$ be an arbitrary sequence of 
$n$-pointed LeBrun metrics, with conformal compactifications 
$(\hat{M}, [\hat{g}_i])$. Assume that there exists a 
constant $\delta_n$ such that
\begin{align}
\label{mest3}
0 < \delta_n \leq Y(\hat{M}, [\hat{g}_i]).
\end{align}
Let $g_{Y,i} \in [\hat{g}_i]$ be a sequence 
of unit volume Yamabe minimizers.
Then there exists a subsequence $g_{Y,j}, \{j\} \subset \{i\},$ 
which converges to either 
$(1)$ a constant scalar curvature metric on a $k$-pointed 
LeBrun orbifold, $1 \leq k \leq n$, or $(2)$ a single 
football metric, that is, $S^4/ \ZZ_m$ with the
round metric, for $2 \leq m \leq n$.  
\end{theorem}
\begin{proof} 

As in Theorem \ref{LBconv}, we choose a compactification with 
fixed conformal factor $\tilde{g}_i \in [ g_i]$, so that the 
sequence $(\hat{M}, \tilde{g}_i)$ will limit to a compactified LeBrun 
orbifold, with each bubble-tree consisting of a string of 
multi-Eguchi-Hanson orbifolds, and 
possibly other AF LeBrun orbifold metrics
(in case some monopole points limit to the boundary of $\H^3$).  
Call this limit space $(\hat{M}_{\infty}, \hat{g}_{\infty})$, and the
finite singular set of convergence $S \subset \hat{M}_{\infty}$.  
We write the sequence of unit volume Yamabe minimizers in each conformal 
class as $g_{Y,i} = v_i^{2} \tilde{g}_i$.  
Let $R_i$ denote the scalar curvature of $g_{Y,i}$. By passing to a
subsequence, assume that $\lim_{i \rightarrow \infty} R_i = R_{\infty}$. 
From the assumption \eqref{mest3}, $R_{\infty} > \delta > 0$. 

The assumption \eqref{mest3} implies that
the Yamabe minimizers  $g_{Y,i}$ satisfy a
uniform Sobolev inequality of the following form   
\cite[Proposition 3.1]{TV2}
\begin{align}
\label{mainSob3}
\Vert f \Vert_{L^{4}} \leq C_S   
\Vert \nabla f \Vert_{L^2} + Vol^{-1/4} \Vert f \Vert_{L^2},
\end{align}
for any $f \in C^{0,1}(\hat{M})$. 
The $L^2$-norm of the curvature is 
uniformly bounded, as seen above in Section \ref{GBHS}. We may 
therefore quote the compactness theorem of \cite[Theorem 1.1]{TV2} 
to obtain a subsequence converging to a 
multi-fold limit. 
The Sobolev inequality \eqref{mainSob3} and $b_1(\hat{M}) = 0$ together 
imply that the limit 
must be an irreducible orbifold \cite[Proposition 7.2]{TV2}. 

Let $C_i = \max_{\hat{M}} u_i$. If $C_i$ remains bounded 
from above, then we argue as follows. 
On any compact subset $D \subset \hat{M}_{\infty} \setminus S$, 
we have a Harnack inequality (since the conformal 
factors are bounded from above, and the sequence is 
smoothly converging away from the singular set). This implies 
that the conformal factors $u_i$ will either have a strictly 
positive limit on $\hat{M}_{\infty}$ or will uniformly crash to 
zero on $\hat{M}_{\infty}$. The latter case cannot happen 
since $Vol(g_{Y,i}) = 1$. 
Consequently, the limit must be a CSC (constant scalar curvature)
metric conformal to 
$(\hat{M}_{\infty}, \hat{g}_{\infty})$, and we are in Case (1). 

So we next assume that $C_i \rightarrow \infty$ as $i \rightarrow \infty$. 
Let $x_i$ be points such that 
$u_i(x_i) = C_i$. 
We first assume that $x_i \rightarrow x \in \hat{M}_{\infty} \setminus S$.
Then the usual conformal dilation argument says that 
a bubble is forming on the smooth part \cite{Schoen3}. 
To summarize, this is by looking at the rescaled functions 
$\tilde{u}_i(y) = C_i^{-1} u_i ( C_i^{-1} y)$.
By elliptic theory, this sequence has a subseqence converging 
to a positive solution of $-6 \Delta u = R_{\infty} u^{3}$ on $\RR^4$,
satisfying $u(0) =1$. From our assumption on the Yamabe invariant, 
$R_{\infty} > \delta > 0$.
By Caffarelli-Gidas-Spruck \cite{CGS}, the limit must be the
spherical metric. This implies 
that the Yamabe invariant satisfies $Y(\hat{M},[\tilde{g}_i]) \gtrsim Y(S^4, [g_S])$,
for $i$ large, which contradicts Theorem \ref{est1} above. 
Therefore, we must have $x_i \rightarrow p 
\in S$ as $i \rightarrow \infty$. 
Let $(X_1,g_1)$ be the first bubble at $p$. 
Assume that, after the rescaling the sequence 
to limit to $(X_1,g_1)$ (pointed convergence based at $p$), 
$x_i$ limits to a finite point of $X_1$. The same 
conformal dilation argument shows that the rescaled conformal 
factor (with  $g_1$ as background metric) must be 
bounded from above. Away from the singular points 
of convergence $p_{1,j'}, j' = 1\dots k_1$, we 
again have a Harnack inequality. So the rescaled 
conformal factor either (a) limits identically to 
zero away from the singlar points, or (b) has a finite 
positive limit everywhere. In Case (b), 
the limit of the original sequence must then be a CSC metric 
on the conformal compactification $(\hat{X}_1, \hat{g}_1)$,
since this is the only possible irreducible 
orbifold limit. If $X_1$ is an Eguchi-Hanson orbifold, this cannot 
happen by Proposition \ref{posprop}. So the only possibility in 
Case (b) is that $(X_1, g_1)$ is an orbifold LeBrun AF metric, 
whose compactification is Case (1). 
Case (a) splits into two possibilities. 
Case (a1) is that $x_i$ will limit to $\infty$ in $(X_1,g_1)$. In this 
case, the only possible irreducible orbifold limit 
will be a metric on the compactification of the 
``neck'' region. This follows since all ALE spaces in the bubble tree 
described in Theorem \ref{LBconv} have one end, the 
only possible neck regions are modeled on $\RR^4 / \ZZ_m$. 
The only CSC metric on the compactification of this 
is the $S^4/\ZZ_m$-football metric, by the 
Obata-Tashiro Theorem \cite{Obata,Tashiro}. Therefore,
Case (a1) is exactly Case (2). 
Case (a2) is that $x_i$ limits to 
one of the singular points of convergence of $X_1$ in 
this scaling. We then repeat the above argument around 
this singular point. In general, we repeat the entire argument at
different scalings to see that the limit must be 
(i) a CSC metric on exactly one of the compactified orbifolds 
in the bubble-tree, or (ii) limit occurring on a neck region. 
For Case (i), the Gibbons-Hawing orbifolds do not 
admit CSC metrics by Proposition \ref{posprop}, so 
Case (i) is exactly Case (1). The argument above 
shows that Case (ii) is exactly Case (2). 

 Finally, in Case (1), the limit can never be a 
compactified $0$-pointed LeBrun metric. 
This is conformal to $S^4$ with the round metric, 
and the only CSC metrics in this conformal class are of constant curvature
\cite{Obata}, so have maximal Yamabe invariant. This would contradict \eqref{uest}.
This also proves $2 \leq m$ in Case (2). 
\end{proof}

\subsection{Completion of proofs of Theorems \ref{n=2thm} and \ref{imaint}}
The Yamabe invariant for $n =2,3$ is stricly positive
by Proposition \ref{23prop}. 
For Theorem \ref{n=2thm}, as  $d_H(p_1, p_2) \rightarrow \infty$,
Theorem \ref{poslim} says the only possible 
limit is the compactified $1$-pointed LeBrun metric,
since Case (2) obviously does not happen. This is conformal 
to the Fubini-Study metric, which is the unique CSC 
metric in its conformal class by Obata's Theorem \cite{Obata}.
We note that the behavior of the Yamabe-minimizer on a connect sum is 
typically non-symmetric \cite{Kobayashi, JoyceCSC}. 
Since the Yamabe minimizers must limit to the Fubini-Study 
metric, it is then obvious that for $d_H(p_1, p_2)$ 
very large, the Yamabe minimizers cannot be
invariant under the conformal involution which flips the 
two monopole points (see Section \ref{sn2}). Therefore, 
there must always be at least two distinct Yamabe minimizers,
which are related by this conformal involution. 

We next address the $K$-symmetric limit to $S^4$. 
The existence of $K$-symmetric minimizers follows from 
\cite{Hebey1996}. Again, by \cite{TV2} we can 
find a subsequence converging to a multi-fold,
but which now may have reducible points. 
But from the estimate on the $K$-Yamabe 
invariant in Theorem \ref{symesti}, 
as $d_H(p_1, p_2) \rightarrow \infty$, it is clear that 
the only possibility for a $K$-Yamabe minimizer is 
$S^4$ with the round metric. This follows because
a reflection interchanging the two monopole point 
is clearly not an isometry of the two Yamabe minimizers,
so these are not possible limits.  From the arguments in the proof of 
Theorem \ref{poslim}, another possible limit is $\CP^2 \vee \CP^2$, 
with the Fubini-Study metric scaled to have $Vol = 1/2$ on each factor.
But this cannot occur since the Yamabe invariant of this limit is $24 \pi
> 8 \pi \sqrt{6}  = Y(S^4, g_S)$, and this 
would contradict Theorem \ref{symesti}. Similarly, 
the limit $\CP^2 \vee S^4$ cannot happen either. 
Therefore, any sequence 
of Yamabe minimizers must concentrate entirely in the neck 
region. The only possible limit is then $S^4$ with the round 
metric, with Burns metrics bubbling off at the 2 singular points
of convergence. This follows since any reducible limit would 
be several $S^4$-s wedged together, which would have Yamabe invariant 
strictly larger than $8 \pi \sqrt{6}$.

\begin{remark}{\em
As shown in \cite{HV}, $G$ is generated by 
$K$ and an extra involution $\Lambda$ which is not a 
lift of a hyperbolic isometry, and is quite difficult 
to describe explicitly.   
It is likely that as $d(p_1, p_2) \rightarrow \infty$, the 
$G$-symmetric metric must also limit to $S^4$.  
But since the test metric in \eqref{smetric} is not invariant 
under $\Lambda$, we cannot say this for certain.} 
\end{remark}

The fourth metric in Theorem \ref{n=2thm} is obtained by 
adapting the CSC-gluing argument of Joyce to this 
problem \cite{JoyceCSC}, see also \cite{MPU}. 
Recall that the Euclidean Schwarzschild metric in dimension $n$
is defined as 
\begin{align}
g = \left(1 + \frac{m}{(n-1)r^{n-2}} \right)^{\frac{4}{n-2}} g_0,
\end{align}
on $\RR^n \setminus \{ 0 \}$, where $g_0$ is the Euclidean metric,
and $m > 0$ is the mass parameter. This metric is scalar-flat, 
locally conformally flat, and AF with two ends. One chooses a 
conformal factor which is close to 
the Fubini-Study metric on neighborhoods of the 
monopole points, and is close to a scaled-down Schwarzschild 
neck region in between; this will be the approximate CSC-metric. 
One then uses the implicit function theorem, together with the crucial 
fact that the Fubini-Study metric is CSC-nondegenerate, to perturb to 
a CSC metric. This adaptation is straightforward, but since 
the argument is quite lengthy, we omit the details due to space considerations. 

As $d_H(p_1, p_2) \rightarrow 0$, Case (1) in Theorem \ref{poslim}
could only be a compactified LeBrun metric with a
single point of multiplicity $2$. This 
is conformal (minus the orbifold point) to the Eguchi-Hanson metric, 
which does not admit any CSC metric by Theorem \ref{nonexist}, 
so this case cannot happen. Consequently, the only possibility for the 
limit is the $S^4 / \ZZ_2$-football with the round metric
(since the limit of Yamabe minimizers must be irreducible). 
Similarly, Theorem \ref{nonexist} implies that as 
$d_H(p_1, p_2) \rightarrow 0$, the only 
possible $K$-symmetric limit is the $S^4/ \ZZ_2$-football 
with the round metric. To see this, again the only concentration 
can occur in the neck region, which is $\RR^4 / \ZZ_2$. 
The limit no longer has to be irreducible. But the lowest 
energy reducible limit would be the wedge of two $S^4/\ZZ_2$-footballs,
whose Yamabe energy is $8 \pi \sqrt{6}$. However, since in this 
case the points are not limiting to the boundary of $\H^3$, 
Theorem \ref{symesti} shows that the Yamabe energy must be 
{\em{strictly}} less than $8 \pi \sqrt{6}$, so this cannot happen, 
and therefore the limit must be the irreducible $S^4/\ZZ_2$-football.

For the first case in Theorem \ref{imaint}, as $p_3 \rightarrow p_2$, 
the only possible limits from Theorem \ref{poslim}
are Case (1): a compactified LeBrun metric with a single 
monopole point, and another multiplicity $2$ point, 
or Case (2): a $S^4/ \ZZ_2$-football with the round metric. 
However, under the assumptions of Theorem \ref{imaint}, 
\eqref{lthmest} says the Yamabe invariant 
is strictly less than that of the $S^4/ \ZZ_2$-football,
so the limit must be Case (1). 
  
 Finally, for the last case in Theorem \ref{imaint}, as
$p_3 \rightarrow \partial \H^3$, the Yamabe invariant is strictly 
less than that of $g_{FS}$ by Theorem \ref{est1}. 
Thus the only possible limit is Case (1), a $2$-pointed 
LeBrun compactified metric. 

\section{Questions}
We conclude with a list of questions. 
\begin{itemize}
\item{
What is the optimal lower bound for the Yamabe invariant 
in \eqref{mest} for $n \geq 2$?
We conjecture that 
$\delta_n = Y_{orb}( S^4 / \ZZ_n, [g_{S}])$; this is true for $n =2$, 
as seen above in Theorem \ref{n=2thm}. 
Furthermore, the assumption that the points 
are contained in a compact set $\mathcal{K} \subset \mathcal{H}^3$
should not be necessary. Removing this assumption would 
imply that the moduli space of Yamabe minimizing LeBrun metrics has a nice 
compactification for any $n$. As seen above, this is true for $n = 2, 3$.}
\vspace{3mm}
\item{Does the compactified LeBrun metric on 
$\mathcal{O}(-n)$ admit a CSC metric for $n \geq 3$? 
The answer is no for $n = 2$ since this is the 
compactified Eguchi-Hanson metric, which was ruled out by Theorem \ref{nonexist}. 
For $n \geq 3$, Theorem~\ref{nonsym} rules out any  
symmetric solution, but is there a non-symmetric solution?
Furthermore, 
as shown in Theorem \ref{LBconv}, if all monopole points approach 
a single point at a uniform rate, the bubble-tree structure is 
a compactified LeBrun negative-mass metric on $\mathcal{O}(-n)$, 
with a Gibbons-Hawking multi-Eguchi-Hanson bubbling 
off. From the arguments in Section \ref{convergence}, the limit 
of the Yamabe minimizers could limit to either a CSC metric on 
compactified $\mathcal{O}(-n)$, or to the $S^4/ \ZZ_n$-football 
(since the compactified GH metric does not admit a CSC metric 
by Theorem \ref{nonexist}). Which one actually happens for $n \geq 3$?}
\vspace{3mm}
\item{Except for the case of a single monopole point with multiplicity 
$n$, does a $[\hat{g}_{LB}]$ orbifold conformal
class always admit a CSC metric? } 
\vspace{3mm}
\item{In Theorem \ref{n=2thm}, for $n=2$ we determined the limiting 
behavior of the $K$-Yamabe minimizers as 
$d(p_1,p_2) \rightarrow \infty$ (the limit is $S^4$). Recall that $K$ is an index 
$2$ subgroup of full conformal group $G$. What is the 
limit of the $G$-Yamabe minimizers? }
\vspace{3mm}
\item{
Are CSC metrics on compactified LeBrun metrics 
CSC nondegenerate? If so, then it would then be possible to 
apply the Joyce gluing technique to obtain more non-Yamabe-minimizing
examples. }
\vspace{3mm}
\item{
LeBrun metrics with torus action are a special 
case of Joyce metrics \cite{Joyce1995}. These depend on a choice 
of points on the boundary of hyperbolic $2$-space. What happens 
as these metrics degenerate?}
\end{itemize}
\providecommand{\bysame}{\leavevmode\hbox to3em{\hrulefill}\thinspace}

\end{document}